\documentstyle[11pt,a4,amscd,amssymb,verbatim]{amsart} 
\newtheorem{Thm}{Theorem}
\newtheorem{Rem}{Remark}
\newtheorem{Cor}{Corollary}  
\newtheorem{lemma}{Lemma}

\newcommand{\Char}{\operatorname{Char}}
\newcommand{\coeff}{\operatorname{coeff_{E_{\bf n}}}}

\newcommand{\Ima}{\operatorname{Im}}

\newcommand{\End}{\operatorname{End}}

\newcommand{\Span}{\operatorname{span}}

\newcommand{\ind}{\operatorname{ind}}
\newcommand{\res}{\operatorname{res}}

\newcommand{\C}{{\mathbb C}}

\newcommand{\A}{{\cal A}}

\edef\savecatcodeat{\the\catcode`@}
\catcode`\@=11

\def\tb@ifSpecChars#1#2{#1}
\def\tb@ifNoSpecChars#1#2{#2}

\def\tableau{%
  \bgroup
  \@ifstar{\let\Tif\tb@ifNoSpecChars\tb@tableauB}
          {\let\Tif\tb@ifSpecChars\tb@tableauB}}

\def\tb@tableauB{
  \@ifnextchar[{\tb@tableauC}{\tb@tableauC[]}}

\def\tb@tableauC[#1]{\hbox\bgroup%
    \let\\=\cr
    \def\bl{\global\let\tbcellF\tb@cellNF}%
    \def\tf{\global\let\tbcellF\tb@cellH}
    \dimen2=\ht\strutbox \advance\dimen2 by\dp\strutbox%
    \ifx\baselinestretch\undefined\relax%
    \else%
       \dimen0=100sp \dimen0=\baselinestretch\dimen0%
       \dimen2=100\dimen2 \divide\dimen2 by\dimen0%
    \fi%
    \let\tpos\tb@vcenter
    \tb@initYoung
    \tb@options#1\eoo
    \let\arrow\tb@arrow%
    \dimen0=\Tscale\dimen2%
    \dimen1=\dimen0 \advance\dimen1 by \tb@fframe%
    \lineskip=0pt\baselineskip=0pt
    \def\tb@nothing{}%
    \def\endcellno{$\rss\egroup\bss\egroup}
    \def\endcell{\endcellno\kern-\dimen0}
    \def\begincell{\vbox to\dimen0\bgroup\vss\hbox to\dimen0\bgroup\hss$}%
    \let\overlay\tb@overlay%
    \let\fl\tb@fl%
    \let\lss\hss\let\rss\hss\let\tss\vss\let\bss\vss
    \def\mkcell##1{
        \let\tbcellF\tb@cellD
        \def\tb@cellarg{##1}
        \ifx\tb@cellarg\tb@nothing\let\tb@cellarg\tb@cellE\fi%
        \begincell\tb@cellarg\endcellno
        \tbcellF}
    \let\savecellF\tbcellF
     \Tif{\catcode`,=4\catcode`|=\active}{}\tb@tableauD}%

\let\tb@savetableauD\tableauD
{
    \catcode`|=\active \catcode`*=\active \catcode`~=\active%
    \catcode`@=\active
\gdef\tableauD#1{%
  \Tif{
    \mathcode`|="8000 \mathcode`*="8000%
    \mathcode`~="8000 \mathcode`@="8000%
    \def@{\bullet}%
    \let|\cr
    \let*\tf
    \let~\sk
  }{}%
  \tpos{\tabskip=0pt\halign{&\mkcell{##}\cr#1\crcr}}%
  \global\let\tbcellF\savecellF
  \egroup
  \egroup}
}
\let\tb@tableauD\tableauD
\let\tableauD\tb@savetableauD
\let\tb@savetableauD\undefined

\def\tb@options#1{\ifx#1\eoo\relax\else\tb@option#1\expandafter\tb@options\fi}

\def\tb@option#1{%
  \if#1t\let\tpos\tb@vtop\fi
  \if#1c\let\tpos\tb@vcenter\fi
  \if#1b\let\tpos\vbox\fi
  \if#1F\tb@initFerrers\fi
  \if#1Y\tb@initYoung\fi
  \if#1s\tb@initSmall\fi
  \if#1m\tb@initMedium\fi
  \if#1l\tb@initLarge\fi
  \if#1p\tb@initPartition\fi
  \if#1a\tb@initArrow\fi
}

\def\tb@vcenter#1{\ifmmode\vcenter{#1}\else$\vcenter{#1}$\fi}

\def\tb@vtop#1{\hbox{\raise\ht\strutbox\hbox{\lower\dimen0\vtop{#1}}}}

\def\tb@initPartition{\def\Tscale{.3}}
\def\tb@initSmall{\def\Tscale{1}}
\def\tb@initMedium{\def\Tscale{2}}
\def\tb@initLarge{\def\Tscale{3}}

\def\tb@initArrow{\dimen2=1.25em}

\def\tb@initYoung{%
  \def\tb@cellE{}
  \let\tb@cellD\tb@cellN
  \def\sk{\global\let\tbcellF\tb@cellNF}}
\def\tb@initFerrers{%
  \def\tb@cellE{\bullet}
  \let\tb@cellD\tb@cellNF
  \def\sk{\bullet}}

\tb@initMedium

\def\tb@sframe#1{%
  \vbox to0pt{
    \vss
    \hbox to0pt{%
      \hss
      \vbox to\dimen1{
        \hrule depth #1 height0pt
        \vss
        \hbox to\dimen1{
          \vrule width #1 height\dimen1
          \hss
          \vrule width #1
          }%
        \vss
        \hrule height #1 depth 0in
        }%
      \kern-\tb@hframe
      }%
    \kern-\tb@hframe}}

\def\tb@hframe{.2pt}\def\tb@fframe{.4pt}\def\tb@bframe{1.2pt}
\def\tb@cellH{\tb@sframe{\tb@bframe}}       
\def\tb@cellNF{}                            
\def\tb@cellN{\tb@sframe{\tb@fframe}}       
\let\tbcellF\tb@cellN                       

\def\tb@rpad{1pt}
\def\tb@lpad{1pt}
\def\tb@tpad{1.8pt}
\def\tb@bpad{1.8pt}

\def\tb@overlay{\endcell\@ifnextchar[{\tb@overlaya}{\begincell}}
\def\tb@overlaya[#1]{\vbox to\dimen0\bgroup%
  \tb@overlayoptions#1\eoo%
  \tss\hbox to\dimen0\bgroup\lss}
\def\tb@overlayoptions#1{\ifx#1\eoo\relax\else\tb@overlayoption#1\expandafter\tb@overlayoptions\fi}

\def\tb@overlayoption#1{
  \if#1t\def\tss{\vskip\tb@tpad}\let\bss\vss\fi
  \if#1c\let\tss\vss\let\bss\vss\fi
  \if#1b\def\bss{\vskip\tb@bpad}\let\tss\vss\fi
  \if#1l\def\lss{\hskip\tb@lpad}\let\rss\hss\fi
  \if#1m\let\lss\hss\let\rss\hss\fi
  \if#1r\def\rss{\hskip\tb@rpad}\let\lss\hss\fi
}

\def\tb@fl{\endcell\begincell\vrule depth 0pt width \dimen0 height \dimen0 \endcell\begincell}

\@ifundefined{diagram}{}{
\def\tb@arrowpad{.5}

\newoptcommand{\tb@arrow}{\@ne}[2]{%
  \endcell
   \begingroup%
   \let\dg@getnodesize\tb@getnodesize
   \dg@USERSIZE=#1\relax%
   \ifnum\dg@USERSIZE<\@ne \dg@USERSIZE=\@ne \fi%
   \dg@parse{#2}%
   \dg@label{\tb@draw{#1}{#2}}}

\def\tb@getnodesize#1#2#3#4#5{\dimen3=\tb@arrowpad\dimen2 #4=\dimen3 #5=\dimen3\relax}
\def\tb@getnodesize#1#2#3#4#5{\ifnum#2=0\ifnum#3=0\tb@getnodesizetail{#4}{#5}\else\tb@getnodesizehead{#4}{#5}\fi\else\tb@getnodesizehead{#4}{#5}\fi}
\def\tb@getnodesizetail#1#2{\dimen3=.5\dimen2 #1=\dimen3 #2=\dimen3}
\def\tb@getnodesizehead#1#2{\dimen3=.5\dimen2 #1=\dimen3 #2=\dimen3}

\def\tb@draw#1#2#3#4{%
        \dg@X=0\dg@Y=0\dg@XGRID=1\dg@YGRID=1\unitlength=.001\dimen0%
        \dg@LBLOFF=\dgLABELOFFSET \divide\dg@LBLOFF\unitlength%
        \dg@drawcalc
        \begincell
        \let\lams@arrow\tb@lams@arrow
        \begin{picture}(0,0)\begingroup\dg@draw{#1}{#2}{#3}{#4}\end{picture}%
        \endcell
        \endgroup
        \begincell}
}

\def\tb@lams@arrow#1#2{%
 \lams@firstx\z@\lams@firsty\z@
 \lams@lastx#1\relax\lams@lasty#2\relax
 \lams@center\z@
 \N@false\E@false\H@false\V@false
 \ifdim\lams@lastx>\z@\E@true\fi
 \ifdim\lams@lastx=\z@\V@true\fi
 \ifdim\lams@lasty>\z@\N@true\fi
 \ifdim\lams@lasty=\z@\H@true\fi
 \NESW@false
 \ifN@\ifE@\NESW@true\fi\else\ifE@\else\NESW@true\fi\fi
 \ifH@\else\ifV@\else
  \lams@slope
  \ifnum\lams@tani>\lams@tanii
   \lams@ht\ten@\p@\lams@wd\ten@\p@
   \multiply\lams@wd\lams@tanii\divide\lams@wd\lams@tani
  \else
   \lams@wd\ten@\p@\lams@ht\ten@\p@
   \divide\lams@ht\lams@tanii\multiply\lams@ht\lams@tani
  \fi
 \fi\fi
 \ifH@  \lams@harrow
 \else\ifV@ \lams@varrow
 \else \lams@darrow
 \fi\fi
}

\catcode`\@=\savecatcodeat
\let\savecatcodeat\undefined

\begin{document}    
\large

\title{On the representation theory of an algebra of braids and ties   } 
\author{Steen Ryom-Hansen}
\address{Instituto de Matem\'atica y F\'isica, Universidad de Talca \\
Chile\\ steen@@inst-mat.utalca.cl }
\footnote{Supported in part by Programa Reticulados y Simetr\'ia and by FONDECYT grants 1051024 and 1090701}

\maketitle

\begin{abstract}
We consider the algebra ${\cal E }_n(u)$ introduced by F. Aicardi and J. Juyumaya as an abstraction of the 
Yokonuma-Hecke algebra. 
We construct a tensor space representation for ${\cal E }_n(u)$ and show that this is faithful. We use it to 
give a basis of ${\cal E }_n(u)$ and to classify its irreducible representations.
\end{abstract}

\section{Introduction}
We initiate in this paper a systematic study of the representation theory of an algebra $  {\cal E }_n(u)$
defined by F. Aicardi and J. Juyumaya. 
Let
$ G $ be a Chevalley group over $ {\mathbb F }_q $ with Borel group $ B $ and
maximal unipotent subgroup $ U $. 
The origin of $  {\cal E }_n(u)$ is in the Yokonuma-Hecke 
algebra $  {\cal Y }_n(u)$, 
which is defined similarly as the Iwahori-Hecke algebra 
but with $ B $ 
replaced by $ U $.
That is, $  {\cal Y }_n(u)$
is the endomorphism algebra of the induced $ G$-module $ \ind_U^G 1 $.
Yokonuma gave in [Y] a presentation of $  {\cal Y }_n(u)$
along the lines of the standard $  T_i \, $-presentation of the Iwahori-Hecke algebra, but 
the introduction of 
$  {\cal E }_n(u)$ 
is more naturally motivated by 
the new  presentation of 
$  {\cal Y }_n(u)$ found by Juyumaya in [J2].
For type $ A_n $, this new 
presentation has generators $ T_i,   \, \, i = 1, \ldots,  n-1 $ 
and $ f_i,  \, \, i = 1, \ldots ,  n $
where the $ f_i $ generate a product of cyclic groups and the $ T_i $ satisfy the usual 
braid relation of type $ A $, but do {\it not} coincide with Yokonuma's $ T_i $-generators. The
quadratic relation takes the form 
$$ T_i^2 = 1 + (u-1) e_{i} (1+T_i) $$ 
for $  e_{i} $ a complicated expression involving $ f_i $ and $ f_{i+1} $.

\medskip
The algebra $  {\cal E }_n(u)$ is obtained by 
leaving out the $ f_i $, but 
declaring the $ e_{i} $ new generators, denoted $E_i $. 
It was introduced by Aicardi and Juyumaya in [AJ].
They showed that $  {\cal E }_n(u)$ is finite dimensional and that it has connections to knot 
theory via the Vasiliev algebra. They also constructed 
a diagram calculus for $  {\cal E }_n(u)$
where the $ T_i $ are represented by braids in the usual sense and the $E_i $ by ties. 
Using results from [CHWX], they moreover showed that $ {\cal E }_n(u) $ can be Yang-Baxterized in the sense of V. Jones, [Jo].

\medskip
In this paper we initiate a systematic study of the representation theory of $  {\cal E }_n(u)$, obtaining a complete classification of its simple 
modules for generic choices of the parameter $ u $. In [AJ], this was achieved only for $ n= 2, 3$. 
An interesting feature of this classification is the construction of  
a tensor space module $ V^{\otimes n } $ for $  {\cal E }_n(u)$. It 
was in part inspired by 
the tensor module for the Ariki-Koike algebra in [ATY] --  see also [RH].
A main property of $ V^{\otimes n} $ is its faithfulness that  
we obtain as a corollary to our theorem \ref{Gbasis} giving
a basis $ G $ for 	
$  {\cal E }_n(u)$. The dimension of $  {\cal E }_n(u)$ 
turns out to be $  B_n n! $ where $ B_n $ is the Bell number, 
i.e. the number of set partitions of $ \{1, 2, \ldots , n \} $. 

\medskip
The appearance of the Bell number is somewhat intriguing and may indicate 
a connection to the partition algebra defined independently by 
P. Martin in [M] and V. Jones in [Jo1], 
but as we indicate in the remarks following 
corollary \ref{faithful}, we do not think at present that the connection can be 
very direct.

\medskip
Given the tensor module, the classification of the irreducible modules follows the principles laid out 
in James's famous monograph on the representation theory of the symmetric group, [Ja].

\medskip
Let us briefly explain the organization of the paper. Section 2 
contains the definition of the algebra $  {\cal E }_n(u)$. In section 3 we start out by giving the construction of the tensor space  
$ V^{\otimes n} $. We then construct the subset $ G \subset   {\cal E }_n(u)$ and show that it generates $ {\cal E }_n(u) $. Finally we show that it maps to 
a linearly independent set in $ \End( V^{\otimes n})  $, thereby obtaining the faithfulness 
of $ V^{\otimes n} $ and the dimension of $ {\cal E }_n(u) $.

\medskip
In section 4 we recall the basic representation theory of the symmetric group and the 
Iwahori-Hecke algebra, and use the previous sections to construct certain simple modules 
for $ {\cal E }_n(u) $ as pullbacks of the simple modules of these. In section 5 we show that $ {\cal E }_n(u) $ is selfdual by constructing a nondegenerate invariant 
form on it. This involves the Moebius function for the usual partial order on set partitions. In 
section 6 we give the classification of the simple modules of $ {\cal E }_n(u) $, to a large extent following the approach of James's book, [Ja]. 
Thus, 
we especially introduce a parametrizing set $ {\cal L}_n $ for the irreducible modules, analogues of the permutations modules and prove James's 
submodule theorem in the setup. The simple modules, the Specht modules, turn out to be a 
combination of the Specht modules for the Hecke algebra and for the symmetric group and 
hence $ {\cal E }_n(u) $ can be seen as a combination of these two. Finally, in the last section 
we raise some questions connected to the results of the paper.

\medskip
It is a great pleasure to thank J. Juyumaya for telling me about $  {\cal E }_n(u) $ and for many useful conversations. Thanks are also due to C. Stroppel and A. Ram for useful discussions 
during the ALT-workshop at the Newton Institute for Mathematical Sciences
and to the 
referee for useful comments that helped improving the presentation of the paper.
Finally, it is a special pleasure to thank V. Jones for useful discussions in Talca during 
his one month stay at the Universidad de Talca. During his visit the city of Talca was 
badly affected by an 
earthquake of magnitude 8.8
on the Richter scale, among the highest ever recorded.

\section{Definition of $  {\cal E }_n(u)$}
In this section we introduce the algebra $ {\cal E }_n(u) $, the main object of our work.
Let $ {\cal A} $ be the principal ideal domain $ \C[u,u^{-1}] $ where
$ u$ is an unspecified variable. 
We first define the algebra $ {\cal E }^{\A}_n(u) $ as the associative unital ${\cal A}$-algebra 
on the generators $ T_1, \, \ldots, T_{n-1} $ and 
$ E_1,\,\ldots,  E_{n-1} $ and 
relations 
$$ 
\begin{array}{ll}
(E1) &  T_i T_j = T_j T_i \, \, \,\,\,\, \, \mbox{ if } \,\, \, |i-j | > 1 \\
(E2) & E_i E_j = E_j E_i \, \,\,\,\, \, \forall \, i,j  \\ 
(E3) & E_i T_j = T_j E_i \, \, \,\,\,\, \, \mbox{ if } \,\, \, |i-j | > 1 \\
(E4) & E_i^2  = E_i \\
(E5) & E_i T_i = T_i E_i \\
(E6) &  T_i T_j T_i = T_j T_i T_j \, \, \,\,\,\, \, \mbox{ if } \,\, \, |i-j | =1 \\
(E7) & E_j  T_i T_j = T_i T_j E_i   \, \, \,\,\,\, \, \mbox{ if } \,\, \, |i-j | =1 \\
(E8) &  E_i E_j T_j = E_i T_j E_i =T_j E_i E_j   \, \, \,\,\,\, \, \mbox{ if } \,\, \, |i-j | =1 \\
(E9) & T_i^2   = 1 + (u-1) E_i (1+T_i) 
\end{array}
$$
It follows from $(E9) $ that $ T_i $ is invertible with inverse
$$ T_i^{-1} = T_i +(u^{-1} -1) E_i (1+T_i) $$
so the presentation of $  {\cal E }_n(u)$ is not efficient, since the generators $ E_i $ for $ i \geq 2 $ can be expressed in terms
of $E_1 $. However, for the sake of readability, we prefer the presentation as it stands. 

\medskip
We then define
$ {\cal E }_n(u) $ as 
$$ {\cal E }_n(u): =  {\cal E }^{\A}_n(u) \otimes_{\A} \C(u) $$ 
where $ \C(u) $ is considered as an $ \A $-module through inclusion.

\medskip
This algebra is our main object of study. 
It was introduced by Aicardi and Juyumaya, in [AJ], although the relation $(E9)$ varies 
slightly from theirs since we have changed $ T_i $ to $ -T_i$. They show, among other things, 
that it is finite dimensional.

\medskip
From $ {\cal E }^{\A}_n(u) $ we can consider the specialization to a fixed value $ u_0$ of $ u$ which we denote $ {\cal E }_n({u_0)} $.
However, we shall in this paper only need the case $u_0=1 $, corresponding to 
$$ {\cal E }_n(1) =  {\cal E }^{\A}_n(u) \otimes_{\A} \C $$ 
where $ \C $ is made into an $ \A $-module by taking $ u $ to $1$.
Letting $ S_n $ denote the symmetric group on $ n $ letters, 
there is a natural algebra homomorphism $ \iota: \C S_n \rightarrow {\cal E }_n(1), \, (i,i+1) \mapsto T_i  $. 
It can be shown to be 
injective, using the results of the paper.


\section{The tensor space}
For the rest of the paper we shall write $ K = \C(u) $. 
Let $ V $ be the $K$-vector space 
$$ V = \Span_{K} \{ \, v_i^{j} \, | \, i,j =  1, 2,  \, \ldots, n \,\} $$
We consider the tensor product $  V^{\otimes 2 } $ and define $ E   \in \End_{K} ( V^{\otimes 2 } ) $ by the rules
$$ 
E (v_{i_1}^{j_1} \otimes v_{i_2}^{j_2} ) = \left\{
\begin{array}{ll}
v_{i_1}^{j_1} \otimes v_{i_2}^{j_2}   &  \mbox{ if  }  j_1 = j_2 \\
0 & \mbox{ otherwise } 
\end{array}
\right. 
$$

\noindent
Furthermore we define $ T   \in \End_{K} ( V^{\otimes 2 } ) $ by the rules
$$
T (v_{i_1}^{j_1} \otimes v_{i_2}^{j_2} ) = \left\{ 
\begin{array}{ll}
v_{i_2}^{j_2}  \otimes  v_{i_1}^{j_1}    & \mbox{ if } j_1 \not=  j_2  \\
u \, v_{i_2}^{j_1}  \otimes  v_{i_1}^{j_2}         & \mbox{ if }  j_1 = j_2, \, i_1 = i_2      
\\
v_{i_2}^{j_2}  \otimes v_{i_i}^{j_1}    & \mbox{ if } j_1 = j_2, \,
i_1 <  i_2  \\
u \, v_{i_2}^{j_2}  \otimes v_{i_1}^{j_1}
+ (u-1)\,   v_{i_1}^{j_1}  \otimes v_{i_2}^{j_2}  
   & \mbox{ if } j_1 = j_2, \,
 i_1 >  i_2
\end{array} \right.
$$

\medskip

We extend these operators to operators $ E_i , T_i $ acting in the tensor space $ V^{\otimes n } $ by letting $ E, T $ act 
in the factors $ (i,i+1) $. In other words, $ E_i $ acts as a projection
in the factors at the positions $ (i,i+1) $ with equal upper index,  
whereas $ T_i $ acts as 
a transposition if the upper indices are different and as a Jimbo matrix for the action of the Iwahori-Hecke algebra in the usual tensor 
space if the upper indices are equal, see [Ji].

\begin{Thm} 
With the above definitions $ V^{\otimes n} $ becomes a module for the algebra $ {\cal E}_n(u)$.
\end{Thm} 

\begin{pf*}{Proof}
We must show that the operators satisfy the defining  relations $ (E1),\,  \ldots, (E9) $. 
Here the relations $ (E1),\,  \ldots, (E5) $ are almost trivially satisfied, since $E_i$ acts as a projection.

\medskip
To prove the braid relation $(E6)$ we may assume that $ n= 3$ and must evaluate both sides of $(E6) $ on the basis
vectors $ v_{i_1}^{j_1} \otimes v_{i_2}^{j_2} \otimes v_{i_3}^{j_3} $ of $ V^{\otimes 3 }$. 
The case where $ j_1, j_2, j_3 $ are distinct corresponds to the symmetric group case and (E6) certainly holds.
Another easy case is $ j_1= j_2 = j_3$, where (E6) holds by Jimbo's classical result, [Ji]. 

\medskip
We are then left with the case $ j_1=j_2 \not= j_3 $ and its permutations. In order to simplify notation, we omit 
the upper indices of the factors of the equal $ j$'s and replace the third $ j $ by a prime, e.g. 
$ v_{i_1}^{j_1} \otimes v_{i_2}^{j_2} \otimes v_{i_3}^{j_3} $ is written 
$ v_{i_1} \otimes v_{i_2} \otimes v_{i_3}^{\prime} $ and so on.

\medskip
We may assume that the lower indices of the unprimed factors are $ 1 $ or $ 2 $ since the action of $ T $ just depends 
on the order. Furthermore
we may assume that the lower index of the primed factor is always $ 1 $ since $ T $ always acts as a transposition between a primed and an 
unprimed factor. This gives $12 $ cases. 
On the other hand, 
the cases where the two unprimed factors have equal lower indices are easy, since both sides of 
$ (E6) $ act through $ u \, \sigma_{13}  $, where $ \sigma_{13} $ is the permutation of the first 
and third factor
of the tensor product.
So we are left with the following 6 cases
$$ \begin{array}{lcr}
v_1 \otimes v_2 \otimes v_1^{\prime} & v_1 \otimes v_1^{\prime} \otimes v_2  &   v_1^{\prime} \otimes   v_1 \otimes v_2  \\
v_2 \otimes v_1 \otimes v_1^{\prime} & v_2 \otimes v_1^{\prime} \otimes v_1  &   v_1^{\prime} \otimes   v_2 \otimes v_1  
\end{array}
$$

Both sides of $ (E6)$ act through $ \sigma_{13} $ on the first three of these subcases whereas
the last three subcases involve each one Hecke-Jimbo action. For instance 
$$ T_1 T_2 T_1 ( v_2 \otimes v_1 \otimes v_1^{\prime}) 
= u \,  v_1^{\prime} \otimes  v_1 \otimes v_2 + (u-1 )\,   v_1^{\prime} \otimes  v_2 \otimes v_1 $$
which is the same as acting with $ T_2 T_1 T_2 $. The other subcases are similar.

\medskip
Let us now verify that $(E7)$ holds for our operators. We may once again assume that $ n= 3 $
and must check $(E7)$ on all basis elements 
$ v_{i_1}^{j_1} \otimes v_{i_2}^{j_2} \otimes v_{i_3}^{j_3} $. Once again, the cases of $ j_1,j_2,j_3 $ all distinct 
or all equal are easy. We then need only consider $ j_1 = j_2 \not= j_3$ and its permutations and can once again 
use the prime/unprime notation as in the verification of (E6). 

\medskip
Let us first verify that $ E_1 T_2 T_1 = T_2 T_1 E_2 $. We first observe that $ E_2 $ acts as the identity on exactly 
those basis vectors that are of the 
form $v_{i_1}^{\prime} \otimes   v_{i_2} \otimes v_{i_3}$. Hence
$$  T_2 T_1 E_2 (v_{i_1}^{\prime} \otimes   v_{i_2} \otimes v_{i_3}) =    v_{i_2} \otimes v_{i_3} \otimes v_{i_1}^{\prime} 
= E_1 T_2 T_1 (v_1^{\prime} \otimes   v_{i_2} \otimes v_{i_3}) $$
The missing basis vectors are of the form 
$  v_{i_1} \otimes v_{i_2}^{\prime}  \otimes v_{i_3}$ or 
$  v_{i_1}   \otimes v_{i_2}  \otimes v_{i_3}^{\prime} $ and are hence killed by $ E_2 $ and therefore $ T_2 T_1 E_2 $. 
But one easily checks that they are also killed by $ E_1 T_2 T_1$. 

\medskip 
The relation $ E_2 T_1 T_2 = T_1 T_2 E_1 $ is verified similarly.

\medskip
Let us then check the relation $(E8)$. Once again we take $ n=3$ and consider the action of $ E_1 E_2 T_2 $,  $ E_1 T_2 E_1 $ and 
$ T_2 E_1 E_2 $ in the basis vector
$ v_{i_1}^{j_1} \otimes v_{i_2}^{j_2} \otimes v_{i_3}^{j_3} $.
If the $ j_1, j_2, j_3 $ are distinct, the action of the three operators is zero, and if 
$ j_1 = j_2 = j_3 $ they all act as $ T_2 $. Hence we may once again assume that exactly two of the $j$'s are equal.

\medskip
But it is easy to check that 
each of the three operators acts as zero 
on all vectors of the form $  v_{i_1}^{\prime} \otimes v_{i_2}  \otimes v_{i_3}$, 
$  v_{i_1} \otimes v_{i_2}^{\prime}  \otimes v_{i_3}$ and 
$  v_{i_1} \otimes v_{i_2}  \otimes v_{i_3}^{\prime}$.
and so
we have proved that 
$ E_1 E_2 T_2 = E_1 T_2 E_1 =
 T_2 E_1 E_2 $.

\medskip
Similarly one proves that 
$ E_2 E_1 T_1 =  E_2 T_1 E_2 = T_1 E_2 E_1 $.

\medskip
Finally we check the relation $(E9) $, which by $ (E5) $ can be transferred into
$$  T_i^2   = 1 + (u-1)  (1+T_i) E_i $$ 
It can be checked taking $ n=2$. We consider vectors of the form 
$ v_{i_1}^{j_1} \otimes v_{i_2}^{j_2} $. If $ j_1 \not= j_2 $ then 
$ E_i $ acts as zero and we are done. And if 
$ j_1 = j_2 $, the relation reduces to the usual Hecke algebra square. The theorem is proved.
        
\end{pf*}

Since the above proof is only a matter of checking  relations, it also works over $ {\cal E }^{\A}_n(u) $ and hence we get 
\begin{Rem}
There is a module structure of $ {\cal E }^{\A}_n(u) $ on $ V^{\otimes n}$.
\end{Rem}

Our next goal is to prove that $ V^{\otimes n} $ is a faithful representation of 
$ {\cal E }_n(u) $. Our strategy for this will be to construct a subset $ G $ of $ {\cal E }^{\A}_n(u) $ 
that generates $ {\cal E }^{\A}_n(u) $ as an $ \A$-module and 
maps to a linearly independent subset of $ \End_{\A}(V^{\otimes n}) $ under the representation. We will then also have 
determined the dimension of $ {\cal E }_n(u) $.

\medskip 
Let us start out by stating the following useful lemma. 
\begin{lemma} {\label{formulas}} The following formulas hold in $ {\cal E}_n(u)$ and $ {\cal E }^{\A}_n(u) $.
$$ \begin{array}{c} (a) \,\,\,\,\,\,\,\,\,\, T_j E_i T_j^{-1} = T_i^{-1} E_j T_i  \,\,\,\, \mbox{  if  } \, |i-j]= 1 \\
 \,(b) \,\,\,\,\,\,\,\, \,\,\,T_i^{-1} T_j  E_i  = E_j T_i^{-1} T_j   \,\,\,\, \mbox{  if  } \, |i-j]= 1 \\
(c) \,\,\,\,\,\,\,\,\,\, T_j E_i T_j^{-1} = T_i E_j T_i^{-1}  \,\,\,\, \mbox{  if  } \, |i-j]= 1 
\end{array}
$$
\end{lemma}
\begin{pf*}{Proof}                                  
The formula $(a)$ is just a reformulation of (E7) whereas the formula $(b)$ follows from 
$$ T_i^{-1} = T_i + (u^{-1} -1) E_i(1+T_i) $$ 
combined with (E7) and (E8). Formula $(c)$ is a variation of $(b)$.
\end{pf*}

For $ 1 \leq i < j \leq n $ we define $E_{ij} $ by $  E_i  $ if $ j=i+1$, and otherwise 
$$ E_{ij} := T_i T_{i+1} \ldots T_{j-2} E_{j-1} T_{j-2}^{-1} \ldots T_{i+1}^{-1} T_{i}^{-1} $$
We shall from now on use the notation $ { \bf n} := \{ 1, 2, \ldots , n \} $. For any nonempty
subset $ I \subset  {\bf n}  $ we extend the definition of 
$ E_{ij} $ to 
$$ E_I : = \prod_{ (i, j) \in  I \times I , \,  i<j } E_{ij} $$        
where by convention $ E_I := 1 $ if $ | I | = 1 $.
We now aim at showing that this product is independent of the order in which it is taken.

\medskip
Let us denote by $ s_{i} $ the transposition $ (i,i+1) $. 
Write $ E_{ \{j,k\}} $ for $ E_{ \min\{j,k \},\max\{j,k\}} $. 
Then we have 
\begin{lemma} {\label{reflection}} We have for all $i,j,k$ that
$$ \begin{array}{c} (a) \,\,\,\,\,\,\,\,\,\,\,\,\,\,\,\, T_i E_{jk} T_i^{-1} = E_{ \{s_i j,s_i k\}} \\
(b) \,\,\,\,\,\,\,\,\,\,\,\,\, \,\,\, T_i^{-1} E_{jk} T_i = E_{ \{ s_i j,s_i k\}}
\end{array}
$$
\end{lemma}
\begin{pf*}{Proof}                                  
Let us prove $ (a)$.
We first consider the case where $ i $ is not any of the numbers $  j-1, j, k-1 $ or $  k $. In that case we must 
show that $ T_i $ and $ E_{j,k} $ commute. For $ i < j-1 $ and $ i > k$ this is clear since $ T_i  $ then commutes 
with all 
of the factors of $ E_{j,k} $. And for $ j < i < k-1 $ one can commute $ T_i $ through $ E_{j,k} $ using $(E6)$ and $ (E3)$.

\medskip
For $ i = j-1 $ the formula follows directly from the definition of $ E_{j,k}$. 
For $ i = k $ we get that $ T_i $ commutes with all 
the $ T_l$ factors of $ E_{j,k}$ and hence it reduces to showing that 
$$ T_k E_{k-1} T_k^{-1} = T_{k-1} E_{k} T_{k-1}^{-1} $$
which is true by formula $(c)$ of lemma \ref{formulas}.
For $ i = k-1 $ the formula follows from the definitions and $ (E7)$. 

\medskip
Finally, we consider the case $ i = j $. To deal with this case, we first rewrite 
$ E_{jk} $, using $(c)$ 
of lemma \ref{formulas} 
repeatedly starting 
with the innermost term, in the form 
\begin{equation}{\label{innermost}} E_{jk} = T_{k-1} T_{k-2} \ldots T_{j+1}  E_{j}  T_{j+1}^{-1} \ldots T_{k-2}^{-1} T_{k-1}^{-1} 
\end{equation}
The formula of the lemma now follows from relation $(E7)$.

\medskip
Formula $(b)$ is proved the same way.
\end{pf*}

With this preparation we obtain the commutativity of the factors involved in $E_I $. We have that
\begin{lemma}{\label{commuting_idempotents}} 
The $ E_{ij} $ are commuting idempotents of $ {\cal E}_n(u)$ and $ {\cal E }^{\A}_n(u) $. 
\end{lemma} 
\begin{pf*}{Proof}                                  
The $E_{ij} $ are obviously idempotents in 
$ {\cal E}_n(u)$ and $ {\cal E }^{\A}_n(u) $
so we just have to prove that they commute. 

\medskip
Thus, given $E_{ij} $ and $ E_{kl}$ we show by induction on $ (j-i)+(l-k) $ that they commute with each other. The induction starts 
for $ (j-i)+(l-k) = 2 $, in which case $E_{ij} = E_i  $ and $ E_{kl}= E_k $, that commute by (E2). 

\medskip
Suppose now $ (j-i)+(l-k) > 2  $ and that $E_{ij}, E_{kl}$ is not a pair of the form $ E_{s-1,s+2}, E_{s,s+1}$ for any $ s$. One checks
now there is an $ r $ such that $ E_{ s_r \{i,j\} }, E_{ s_r \{k,l\} } $ is covered by the induction hypothesis. But then, using 
$(a)$ from the previous lemma together with the induction hypothesis, we find that
$$ \begin{array}{c}
E_{ij} E_{kl}  = T_r^{-1}  E_{ s_r \{i,j\} } T_r T_r^{-1}  E_{ s_r \{k,l\} } T_r =
T_r^{-1}  E_{ s_r \{i,j\} }   E_{ s_r \{k,l\} } T_r =  \\
T_r^{-1}    E_{ s_r \{k,l\} }  E_{ s_r \{i,j\} }  T_r =  T_r^{-1}    E_{ s_r \{k,l\} } T_r T_r^{-1}  E_{ s_r \{i,j\} }  T_r = E_{kl} E_{ij} 
\end{array}
$$
as needed. Finally, if our pair is of the form $ E_{s-1,s+2}, E_{s,s+1}$ we use (E8) to finish the proof 
the lemma as follows
$$ 
\begin{array}{l}
E_{s-1,s+2} E_{s,s+1} = T_{s-1} T_{s} E_{s+1} T_{s}^{-1} T_{s-1}^{-1} E_s= 
E_s  T_{s-1} T_{s} E_{s+1} T_{s}^{-1} T_{s-1}^{-1}= \\
E_{s,s+1} E_{s-1,s+2}   
\end{array}
$$ \end{pf*}
We have now proved that the product involved in $ E_I $ is independent of the order taken. We then go on to show that many of the factors 
of this product can be left out.
\begin{lemma} 
Let $ I \subset {\bf n} $ with $ | I | \geq 2 $ and set $ i_0 := \min I $. Then
$$ E_I = \prod_{i: \, i  \in I \setminus \{i_0 \} } E_{i_0 i} $$
\end{lemma} 
\begin{pf*}{Proof}                                  
It is enough to show the lemma for $ I $ of cardinality three. By a direct calculation using the definition 
of $ E_{kl} $ one sees that this case 
reduces to $ I = \{1,2,i\} $. Set now
$$ 
\begin{array}{c}
E^{1}:= E_1 T_1 T_2 \ldots T_{i-1} E_i  T_{i-1}^{-1} \ldots T_{2}^{-1} T_{1}^{-1}  \\
E^{2}:=  T_2 T_3 \ldots T_{i-1} E_i  T_{i-1}^{-1} \ldots T_{3}^{-1} T_{2}^{-1} \,\,\,\,\,\,\,  
\end{array}
$$ 
Then the left hand side of the lemma is $ E^{1} E^{2}$ while the right hand side is $ E^{1} $, so we must show that 
$ E^{1} E^{2} = E^{1} $. But using formula $(a)$ of lemma \ref{formulas} repeatedly this identity reduces to 
$$ E_1 T_1 E_2 T_1^{-1} E_2 =  E_1 T_1 E_2 T_1^{-1} $$ 
which is true by relations (E5) and (E8).
\end{pf*}

In order to generalize the previous results appropriately we need to recall some notation.
A {\it set partition} $ A = \{ I_1, I_2,  \ldots , I_k \}$ 
of $ {\bf n } $ is by definition an equivalence relation 
on $ {\bf n } $ with classes $ I_j$. This means that the 
$ I_j $ are disjoint, nonempty subsets of $ {\bf n } $
with union $ {\bf n } $. We also refer to the $ I_j $ as the {\it blocks} of $ A $. 
The number of distinct 
set partitions of $ {\bf n } $ is called the $n$th Bell number and is written $ B_n $. For example $ B_1 = 1 $, 
$ B_2 = 2 $ and $ B_3 = 5 $. The five set partitions of $ {\bf 3 }=  \{ 1,2, 3 \} $ are 
$$ \{ \{1\},  \{2\},  \{3\} \}, \, \,  \{ \{1\},  \{2, 3\} \},\, \, \{ \{2\},  
\{1, 3\} \}, \, \, \{ \{3\},  \{1,2\} \},
\, \, \{ \{1,2 ,3\} \} $$
Let us denote by $ {\cal P}_n $ the set of all set partitions of $ {\bf n} $. 
There is natural partial order on $ {\cal P}_n $, denoted $ \subset $. It is defined by 
$ A= \{ I_1, I_2,  \ldots , I_k \}  \subset B= \{ J_1, J_2,  \ldots , J_l \} $ if and only if
each $ J_i $ is a union of certain $ I_i $. 

\medskip

Let $ R $ be a subset of $ {\bf n } \times {\bf n} $. Write $ i \smallsmile_R j $ 
if $ (i,j) \in R $ and write
$ \sim_R $ for the equivalence relation induced by $ i \smallsmile_R j $. 
Then $ i \sim_R j $ iff 
$ i = j $ 
or there is a chain $ i = i_1, \,  i_2 , \ldots , i_k = j $ such 
that $ i_s \smallsmile_R  i_{s+1} $ or $ i_{s+1} \smallsmile_R  i_{s} $ for all $ s $. 
Let $ \langle R  \rangle $ denote the set partition corresponding to $ \sim_R $.
For example, if $ R = \emptyset $ we get that $ \langle R  \rangle $ is the trivial 
set partition whose blocks are all of cardinality one.

\medskip
For $ A= \{ I_1, \ldots, I_k \}  \in {\cal P}_n$ we define 
$$ E_A := \prod_i E_{I_i} $$
It follows from lemma {\ref{commuting_idempotents} that the product is independent of 
the order in which it is taken. 

\medskip
For $ w \in S_n $
we 
define $ w A:=  \{ w I_1, w I_2, \ldots, w I_k \} \in {\cal P}_n $. 
If $ w = s_{i_1} s_{i_2} \ldots s_{i_n}   $ is a reduced form of $w$, we write as usual 
$ T_w = T_{i_1} T_{i_2} \ldots T_{i_n} $.
Then we have 
\begin{Cor}{\label{setreflection}} 
With $ A \in {\cal P}_n $ and $ w $ as above the following formula holds
$$ T_w E_A T_w^{-1} = E_{ w  A }$$
\end{Cor} 
\begin{pf*}{Proof} 
This is a consequence  of lemma \ref{reflection} $(a)$ and the definitions.
\end{pf*}

The next lemma is an important ingredient in the construction 
of the basis for $ {\cal E}_n(u)$. 

\begin{lemma}{\label{extension}}
Suppose $ R \subset {\bf n } \times {\bf n} $. 
Then the following formula is valid 
$$ \prod_{i,j: (i,j) \in R } E_{ \{i, j \}} = E_{ \langle R \rangle } $$
\end{lemma} 
\begin{pf*}{Proof} 
Writing 
$ E_R := \prod_{i,j: (i,j) \in R } E_{ \{i,j \}} $ we must 
prove that $ E_R : =E_{ \langle R \rangle } $.
Clearly, all the factors of $ E_R $ are also factors of $ E_{ \langle R \rangle } $. We show that 
the extra factors of 
$  E_{ \langle R \rangle } $
do not change the product of $ E_R $. For this, 
suppose first that the following equations hold for $ i < j < k $
\begin{equation}{\label{equivalence_rel}}
E_{ij} E_{ik} = E_{ij} E_{j k} = E_{ik} E_{j k} = 
E_{ij} E_{j k} E_{i k } 
\end{equation}
Assume now that $ i, j  \in { \bf n } $ satisfy $ i \sim_R j $. Then, by definition, 
there is a chain $ i = i_1, \,  i_2 , \ldots , i_k = j $ with 
$ (i_s,   i_{s+1}) \in R $ or $ (i_{s+1},  i_{s}) \in R $ for all $ s $.
Let $ 1 \leq l < k $ and assume recursively that 
we have $ E_R = E_R E_{ \{  i, i_l    \} }$.  
Then using ({\ref{equivalence_rel}}) we get that also $ E_R := E_R E_{ \{ i,  i_{l+1}  \} } $.
Continuing, we find that 
$ E_R := E_R E_{ \{ i j   \} } $, and so indeed
the extra factors of $  E_{ \langle R \rangle } $
do not change the product $ E_R $.  Thus we are reduced to proving 
({\ref{equivalence_rel}}). 

\medskip
The equation 
$ E_{ij} E_{ik} = E_{ij} E_{j k} E_{i k }  $ was shown in the previous lemma so we only need show that 
$ E_{ik} E_{j k} = E_{ij} E_{j k} E_{i k } $ and $ E_{ij} E_{j k} =  E_{ij} E_{j k} E_{i k } $.

\medskip
We consider the involution $ inv $ of $ {\cal E }^{\A}_n(u) $ given by the formulas
$$ inv(T_i) = T_{n-i} \, \, \, \, \, \, \, \, \, \, inv(E_i) = E_{n-i} $$
Using equation (\ref{innermost}) we find that 
$$ inv(E_{ij}) = E_{n-j, n-i} $$ 
But then $ E_{ik} E_{j k} = E_{ij} E_{j k} E_{i k } $ follows from $ E_{ij} E_{ik} = E_{ij} E_{j k} E_{i k }  $. 

\medskip 
We then show that $ E_{ij} E_{j k} =  E_{ij} E_{j k} E_{i k } $. By the above, it can be reduced to showing 
the identity 
$$  E_{ij} E_{j k} =  E_{ij} E_{i k} $$
Using the definition of the $ E_{ij}$ it can be reduced to the case $ i =1 $, $ j= 2 $, i.e.
$  E_1 E_{2 k} =  E_1 E_{1 k} $.
Using 
formula $(a)$ of lemma \ref{formulas} it becomes the valid identity
$  E_1 E_{2 } =  E_1 T_1 E_{2} T_1^{-1}$, 
\end{pf*}
From the lemma we get the following compatibility between the partial order on $ {\cal P}_n $ 
and the $ E_A$. 
\begin{Cor}{\label{inclusion}} Assume $ A, B \in {\cal P}_n $ 
and let $ C \in {\cal P}_n $ be minimal with respect to $ A \subseteq C $ and $ B \subseteq C $. 
Then $ E_A E_B = E_C $.
\end{Cor}

\medskip
We are now in position to construct the subset $G $ of $ {\cal E }^{\A}_n(u) $.
We define 
\begin{equation}
{\label{basis}}  G:= \{ \, E_A T_w \, | \, A \in {\cal P}_n  , \,  w \in S_n  \} 
\end{equation}
With the theory developed so far we can state the following theorem. 
\begin{Thm} 
The set $ G $ generates $ {\cal E }^{\A}_n(u) $ over $\A$.
\end{Thm} 
\begin{pf*}{Proof}
Consider a word $ w = X_{i_1} X_{i_2} \cdots X_{i_k} $ in the generators $ T_i $ and $ E_i $, i.e. 
$ X_{i_j} = T_{i_j} $ or $ X_{i_j} = E_{i_j} $ for all $ j$. Using 
lemma \ref{reflection}
we can move all the 
$ E_{i} $ to the front position, at each step changing the 
index set by its image under some reflection, and are finally left 
with a word in the $T_{i}$, which is possibly not reduced. If it is not so, it is 
equivalent under the braid relations $(E6) $ to a word with two consecutive 
$ T_i $, see [H] chapter 8. Expanding
the $ T_i^2 $ gives rise to a linear 
combination of $ 1, E_i $ and $ T_i E_i $, where the $ E_i$ can be 
commuted to the front position the same way as before. 
Continuing this way we eventually reach a word in reduced form, that is a 
linear combination of elements of the form 
$ \prod_{(i,j) \in R, \, w \in S_n  } E_{ij} \,  T_w  $
for some subset $ R $ of $ {\bf n } \times {\bf n } $, satisfying $ (i,j ) \in R $ only if 
$ i < j $. 
Using lemma \ref{extension} we may rewrite it as a linear combination of 
$  E_{ \langle R \rangle }  T_w  $ and the proof is 
finished.
\end{pf*}

With these results at hand we can prove the following main theorem.
\begin{Thm}{\label{Gbasis}} 
The set $ G $ is a basis 
of $ {\cal E }^{\A}_n(u) $ and induces bases of 
$ {\cal E}_n(u)$ and $ {\cal E}_n(1)$.
\end{Thm} 
\begin{pf*}{Proof}
By the previous theorem it is enough to show that $ G $ is an $\A$-linearly independent 
subset of $ {\cal E }^{\A}_n(u) $ and induces $K $ and $\C $-linearly independent subsets
of $ {\cal E}_n(u)$ and $ {\cal E}_n(1)$.

\medskip
Assume that there exists a nontrivial linear dependence 
$ \sum_{g \in G } \lambda_g \, g = 0 $ 
where $ \lambda_g \in \A $ for all $ g $.
Let $ \lambda \in \A $ be the greatest common divisor of the $ \lambda_g $ 
and write $ \lambda = (v-1)^M \lambda_1 $ with $ \lambda_1 \in \A $ and $ \lambda_1(1) \not= 0 $. 
Setting $ \mu_g := \lambda_g / (v-1)^M  \in \A $ we obtain an $\A $-linear dependence
$ \sum_{g \in G } \mu_g \, g = 0 $ satisfying $ \mu_g (1) \not= 0 $ for at least one $ g $.
By specializing, we obtain from this a nontrivial $\C $-linear dependence
$ \sum_{g \in G } \mu_g(1) \, g = 0 $ in $ {\cal E }_n(1) $.

\medskip
Denoting by 
$ \psi:  {\cal E }^{\A}_n(u) \rightarrow  \End_{\A}( V^{\otimes n}) $ the representation 
homomorphism we get by specializing a 
homomorphism
$ \psi_1:  {\cal E }_n(1) \rightarrow  \End_{\C}( V^{\otimes n}) $.
We use it to obtain the nontrivial linear dependence
$ \sum_{g \in G } \mu_g(1) \, \psi_1 (g) = 0 $ in $ \End_{\C}( V^{\otimes n}) $.
It is now enough to show that $  \{\, \psi_1 (g) \, | \, g \in G \,\} $ 
is a $ \C$-linearly independent set 
of $ \End_{\C}( V^{\otimes n}) $.

\medskip
But for $u=1$, the action of $ T_i$ in $ V^{\otimes n } $ is just permutation 
of the factors $ (i,i+1)$. Hence, in this case, 
$ E_{kl}$ acts as a projection in the space of equal upper indices in the $kl$'th factors of $ V^{\otimes n } $. In formulas
$$
\begin{array}{l}
 E_{kl} (v_{i_1}^{j_1} \otimes \ldots \otimes v_{i_k}^{j_k} \otimes \ldots \otimes  v_{i_l}^{j_l} \otimes \ldots \otimes v_{i_n}^{j_n})=\\
\left\{
\begin{array}{ll}
 v_{i_1}^{j_1} \otimes \ldots \otimes v_{i_k}^{j_k} \otimes \ldots \otimes  v_{i_l}^{j_l} \otimes \ldots \otimes v_{i_n}^{j_n} & 
\mbox{ if } j_k = j_l \\
0 & \mbox{ otherwise } 
\end{array}
\right.
\end{array}
$$
Thus, for a set partition $ A = \{I_1, I_2 \ldots ,I_s \} \in {\cal P}_n $ 
we get that $E_A $ acts 
as the projection $ \pi_A $ on the space of equal upper indices in 
factors corresponding to each of the $ I_k$. In formulas 
$$
\begin{array}{l}
 E_A (v_{i_1}^{j_1} \otimes \ldots \otimes v_{i_r}^{j_r} \otimes \ldots \otimes  v_{i_s}^{j_s} \otimes \ldots \otimes v_{i_n}^{j_n})=\\
\left\{
\begin{array}{ll}
0  \mbox{ if there exist }  {r,s,k} \mbox{ such that } r,s \in I_k \, \mbox{ and } j_r \not=j_s &  \\
 v_{i_1}^{j_1} \otimes \ldots \otimes v_{i_r}^{j_r} \otimes \ldots \otimes  v_{i_s}^{j_s} \otimes \ldots \otimes v_{i_n}^{j_n} 
\mbox{ otherwise } &
\end{array}
\right.
\end{array}
$$
Let us now consider a linear dependence:
\begin{equation}{\label{linear}}
 \sum_{w \in S_n , \,A \in {\cal P }_n  }\lambda_{w,A} \,T_w \pi_A  = 0 
\end{equation}
with $ \lambda_{w,A} \in \C $.
Take $ A_0 \in {\cal P}_n $ such that $ \lambda_{w,A_0 } \not=0 $ for some 
$ w \in S_n  $ and $A_0 $ is minimal with respect to this condition, where 
minimality refers to the partial order on $ {\cal P}_n $ introduced above.
Suppose that 
$ A_0 = \{I_1, I_2, \ldots ,I_s \} $.
If we take a basis vector of $ V^{ \otimes n } $ 
$$v^{A_0}=v_{i_1}^{j_1} \otimes \ldots \otimes v_{i_k}^{j_k} \otimes \ldots \otimes  v_{i_l}^{j_l} \otimes \ldots \otimes v_{i_n}^{j_n}$$
such that $ j_k= j_l $ if and only if $ k,l $ belong to the same $ I_i$, then 
we get on evaluation in (\ref{linear}), using the minimality
of $ A_0 $, that 
$$
 \sum_{w \in S_n  } \lambda_{w,A_0} \,T_w v^{A_0}  = 0 
$$
We now furthermore take $v^{A_0}$ such that its lower $i$-indices are all distinct. But then 
$ \{ T_w v^{A_0} , w \in S_n  \} $ is a linearly independent set and we conclude 
that $ \lambda_{w,A_0 } =0 $ for all $ w$,
which contradicts the choice of $ A_0$.

\medskip
This shows that the set $ \{ T_w \pi_A \, | \, w \in S_n  , A \in {\cal P }_n \} $ is linearly 
independent. To get
the linear independence of $ \{  \pi_A  T_w \, | \, w \in S_n  , A \in {\cal P }_n \} $ 
we apply corollary \ref{setreflection}.

\medskip
We have shown that $ G $ induces a $ \cal \C $-independent subset of 
$ {\cal E }_n(1) $ and we then conclude, as described above, that it is an 
$  \A $-independent subset of 
$ {\cal E }_n^{\A}(u) $. Since $ K $ is the quotient field of $ \cal A $ it also 
induces a $ K $-independent subset of 
$ {\cal E }_n(u) $ and the theorem is proved.
\end{pf*}

\begin{Cor}
We have $\dim {\cal E }_n(u) = n! B_n$, where $ B_n $ is the Bell number, i.e. the number of set partitions of $ \bf n $.
For example $\dim {\cal E }_2(u) = 4 $, $\dim {\cal E }_3(u) = 30$, etc.
\end{Cor}

\medskip 
The appearance of set partitions in the above, notably corollary {\ref{inclusion},
might indicate a connection between $ {\cal E }_n(u)  $ 
and the partition algebra $ A_n (K) $ introduced independently by 
P. Martin in [M] and V. Jones in [Jo1], see also [HR] for an account of the 
representation theory of $ A_n (K) $. On the other hand,    
the special relation $ (E9 ) $ of $ {\cal E }_n(u)  $ does 
complicate the direct comparison 
$ {\cal E }_n(u)  $ with 
known variations of the partition algebra
and at present we do not believe 
that there can be any straightforward connection. 
The relation $ (E9 ) $ 
reveals the origin of $ {\cal E }_n(u)  $ in the Yokonuma-Hecke algebra.  
Since $ u \not= 1 $, it behaves like a kind of skein relation
in the diagram calculus for $ {\cal E }_n(u)  $, 
which seems awkard to interpret in a partition algebra context. 
Note that $ {\cal E }_n(u)  $ becomes infinite dimensional if $ (E9) $ is left out.

\begin{Cor}{\label{faithful}}
The tensor space $ V^{\otimes n} $ is a faithful $ {\cal E }_n(u) $-module.
\end{Cor}
\begin{pf*}{Proof}
We proved that $ G $ is a basis of $ {\cal E }_n(u)  $ that maps to a linearly independent set in $ \End_K (V^{\otimes n}) $.
\end{pf*}

\medskip

\section{Representation theory, first steps}
We initiate in this section the representation theory of $ {\cal E }_n(u)$. We construct two 
families of irreducible representations 
of $ {\cal E }_n(u)$ as pullbacks of irreducible representations of 
the symmetric group and of the Hecke algebra. 

\medskip
Let $ I \subset  {\cal E }_n(u)$ be the two-sided ideal generated by $ E_i $ for all $ i$; actually $ E_1$ is enough to generate $ I $.
Let furthermore $ J \subset  {\cal E }_n(u)$ be the two-sided ideal generated by $ E_i -1  $ for all $ i $; once again $ E_1 -1 $ is enough to generate $ J $. Recall that $ S_n  $ denotes the
symmetric group on $ n $ letters. 
Let $ H_n(u) $ be the Hecke algebra over $ K $ of type $ A_{n-1} $. It is the $ K$-algebra
generated by $ T_1, \ldots , T_{n-1} $ with relations
$  T_i T_j = T_j T_i \mbox{ if }  | i-j | > 1 $ and 
$$ T_i T_{i \pm 1 } T_i = T_{i \pm 1} T_i T_{i \pm 1 }, \,\,\, \,\, \,   (T_i-u)(T_i+1) =0  $$ 
where $ i $ is any index such that the expressions make sense.
\begin{lemma} 
a) There is an isomorphism $\varphi:  K S_n   \rightarrow {\cal E }_n(u)/I , \, 
s_i \mapsto T_i $.  \newline
b) There is an isomorphism $\psi:  H_n(u)  \rightarrow {\cal E }_n(u)/J, \, 
T_i \mapsto T_i $.
\end{lemma}
\begin{pf*}{Proof}                                  
We first prove $a)$. In $ {\cal E }_n(u) /I $ we have $ T_i^2 =1 $ and hence we obtain  
a surjection $ \varphi:  K S_n   \rightarrow {\cal E }_n(u)/I $ by mapping 
$ s_i $ to $ T_i $.
Consider once again the vector space $ V = \Span_{K} \{ v_i^{j} \, | \, i,j =  1, \,  \ldots , n \} $
and its tensor space $ V^{\otimes n } $ as a representation of $ {\cal E }_n(u)$.
We consider the following subspace $ M \subset  V^{\otimes n } $. 
$$ M = 
\Span_{K} \{ \, v_{i_1}^{j_1} \otimes \ldots \otimes  v_{i_n}^{j_n} \, | \mbox{ the upper indices are all distinct }\} $$
It is easy to check from the rules of the action of ${\cal E }_n(u) $
that $ M $ is a submodule of $ V^{\otimes n } $.
Since the $ E_i $ act as zero in $ M $ we get an induced homomorphism $ \rho: {\cal E }_n(u)/I \rightarrow \End_K(M) $,
where $ \rho(T_i) $ is the switching of the $ i $'th and $ i+1 $'th factors of the tensor product. But then the image of 
$ \rho \circ \varphi $ 
has dimension $ n! $ and we conclude that $ \varphi $ indeed is an isomorphism.

\medskip
In order to prove $b)$ we basically proceed in the same way.
In the quotient ${\cal E }_n(u)/J $ we have 
$ T_i^2 = 1 + (u-1)(1+T_i ) $ which implies the existence of a surjection 
$ \psi:  H_n(u)  \rightarrow {\cal E }_n(u)/J $ mapping $ T_i $ to $ T_i$.	
To show that 
$ \psi $ is injective we this time consider the submodule 
$$ N = 
\Span_{K} \{ \, v_{i_1}^{j_1} \otimes \ldots \otimes v_{i_n}^{j_n} \, | \mbox{ the upper indices are all equal to } 1\} $$
All $ E_i $ act as $ 1 $ in $ N $ and so we get a induced map $ \rho^{\prime} : {\cal E }_n(u)/J \rightarrow \End_K(N) $.
The composition $ \rho^{\prime} \circ \psi $ is the regular representation of $ H_n(u) $ and hence $ \dim \Ima 
(\rho^{\prime} \circ \psi ) = n! $ which proves that also $ \psi $ is an isomorphism.

\end{pf*}
We now recall the well known basic representation theory of $ K S_n  $ and of $ H_n(u) $. 
Let $ \lambda = (\lambda_1, \lambda_2 , \ldots , \lambda_k ) $ be an integer partition of $ 
| \lambda | :=n $ 
and let $ Y(\lambda) $ be its Young diagram. Let $ t^{\lambda} $ (resp. $ t_{\lambda} $) be the $\lambda$-tableau in which the numbers 
$ \{1, 2, \ldots,n \} $ are filled in by rows (resp. columns).
Denote by $ R(\lambda) $ (resp. $ C(\lambda) $) the row (resp. column) stabilizer
of $ t^{\lambda}$.
Define now 
$$ r_{\lambda} = \sum_{ w \in R(\lambda) }  w, \,\,\,\,\,\,\,\,\,\,\,
c_{\lambda} = \sum_{ w \in C(\lambda) }  (-1)^{ l(w) }   w, \,\,\,\,\,\,\,\,\,\,\,
s_{\lambda} = c_{\lambda} r_{\lambda} $$
Then $ s_{\lambda} $ is the Young symmetrizer and $ S(\lambda) = K S_n  s_{\lambda} $ is the Specht module associated with
$ \lambda $. Since $ \Char K = 0 $, the Specht modules are simple and classify the simple modules of $ K S_n $.

\medskip
To give the Specht modules for $ H_n(u) $, we use Gyoja's Hecke 
algebra analogue of the Young symmetrizer, [G], [Mu]. In our setup it looks as follows:
For $ X \subset S_n $, define 
$$ \iota(X) = \sum_{w \in X }  T_w, \, \, \, \, \, \, \, \, \,
\epsilon(X) = \sum_{w \in X } (-u)^{-l(w)} T_w \, \, $$
If for example $ X = S_n  $, we have
$$ T_w \, \iota(S_n ) = u^{ l(w)}  \iota(S_n), \,\,\,\,\,\,\,\, 
T_w \, \epsilon(S_n) = (-1)^{ l(w)}\epsilon(S_n)  $$
for all $ T_w $. We now define
$$ x_{\lambda} = \iota(R(\lambda)), \,\,\, \, \, \, \,\,\,\,\,\,\, \, \, \ 
y_{\lambda} = \epsilon(R(\lambda)) $$
Let $w_{\lambda} \in S_n $ be the element such that $ w_{\lambda} \, t^{\lambda} = t_{\lambda} $.
Then the Hecke algebra analogue of the Young symmetrizer is 
$$ e_{\lambda} = T_{ w_{\lambda}^{-1} } y_{\lambda^{\prime} } T_{ w_{\lambda} } x_{\lambda} = c_{\lambda}(u)   r_{\lambda}(u) $$
where $ c_{\lambda}(u) := T_{ w_{\lambda}^{-1} } y_{\lambda^{\prime} } T_{ w_{\lambda} }  $ and 
$ r_{\lambda}(u) := x_{\lambda}(u) $. 
The permutation module and the Specht module associated with $ \lambda $ are 
defined as $ M_u(\lambda ) := H_n(u) x_{\lambda} $ and $ S_u(\lambda) = H_n(u) e_{\lambda}$.
Since $ u $ is generic, $ S_u(\lambda ) $ is irreducible.

\medskip
For future reference, we recall the following result, see eg. [DJ], [Mu].
\begin{lemma}{\label{futurereference}}
Suppose that $ c_{\lambda}(u) M_u(\mu ) \not= 0 $. Then $ \mu \trianglelefteq \lambda  $.
\end{lemma}
Here $ \trianglelefteq $ refers to the dominance order on partitions of $ n$, 
defined by $ \lambda = (\lambda_1, \lambda_2, \ldots )   \trianglelefteq  \mu = (\mu_1, \mu_2 , \ldots ) $ iff 
$ \lambda_1 + \lambda_2 + \ldots + \lambda_i \le  \mu_1 + \mu_2 + \ldots + \mu_i $  for all $ i $. The dominance order is only a partial order, but we shall 
embed it into the total 
order $ < $ on partitions of $ n $, defined by 
$ \lambda = (\lambda_1, \lambda_2, \ldots )  <  \mu = (\mu_1, \mu_2 , \ldots ) $ iff 
$ \lambda_1 + \lambda_2 + \ldots + \lambda_i \leq  \mu_1 + \nu_2 + \ldots + \mu_i $ for some $ i $ and 
$ \lambda_1 + \lambda_2 + \ldots + \lambda_j = \mu_1 + \mu_2 + \ldots + \mu_j $ for $ j < i $.
We extend $ < $ to a total order on all partitions by declaring $ \lambda < \mu $ if $ | \lambda | < | \mu | .$

\medskip
It is known that 
$ y_{ \lambda^{\prime} } T_w x_{ \lambda} \not= 0 \, \mbox{ only if } w= w_{\lambda} $
see [DJ], [Mu]. Using it we find that 
\begin{equation}{\label{schur_hecke}}
c_{ \lambda}(u) z r_{ \lambda}(u)  = C_z c_{ \lambda}(u)  r_{ \lambda}(u) 
\mbox{ for all } z \in H_n(u) 
\end{equation}
for a constant $ C_z \in K$. It follows that $ s_{\lambda}(u) $ is a preidempotent, 
i.e. an idempotent up to a nonzero scalar. 
There is a similar formula 
\begin{equation}{\label{schur_sym}}
c_{ \lambda} z r_{ \lambda}  = C_z c_{ \lambda} r_{ \lambda} 
\mbox{ for all } z \in K S_n
\end{equation}
in the symmetric group case.

\medskip 
Using the Specht module $ S(\lambda) $ for $ K S_n $ or 
$ S_u(\lambda) $ for 
$ H_n(u) $ we use $ \varphi $ 
or $ \psi $ to obtain 
a simple module for $ {\cal E}_n(u) $, by pulling back. 
On the other hand, these two series of simple modules 
do not exhaust all the simple modules for $ {\cal E}_n(u) $ as we shall see in the next sections.

\section{ ${\cal E}_n(u)^{\prime}$ as a ${\cal E}_n(u)$-module } 
In this section we return to ${\cal E}_n(u)$. We show that it is selfdual as a left module over ${\cal E}_n(u)$ itself.
As a consequence of this we get 
that all simple modules occur as left ideals in ${\cal E}_n(u)$.

\medskip
Denote by $ \ast: {\cal E}_n(u) \rightarrow {\cal E}_n(u)$ the $ K $-linear antiautomorphism
given by $ T_i^{\ast} = T_i $ and
$ E_i^{\ast} = E_i $. To check that $ \ast $ exists we must verify that $ \ast $ leaves the defining relations 
$ (E1),\ldots , (E9) $ invariant. This is obvious for all of them, except possibly for $ (E7) $ where 
it follows by interchanging $ i$ and $ j$. There is a similar antiautomorphism for 
$ {\cal E}_n(1) $, also denoted $ \ast $.

\medskip
We now make the linear dual $ {\cal E}_n (u)^{\prime}$ of $ {\cal E}_n (u) $ into a left $ {\cal E}_n (u)$-module using $ \ast $:
$$ (x f)(y) := f( x^{\ast} y  ) \,\,\,\,\,\,\,\, \mbox{ for } x, y \in  {\cal E}_n (u) , f \in  {\cal E}_n (u)^{\prime}$$

We  need to consider the linear map $$ \epsilon: {\cal E}_n(u) \rightarrow K, \,\,x \mapsto \coeff (x) $$ where
$ \coeff (x) $ is the coefficient of $ E_{ \bf n}$ when 
$ x \in {\cal E}_n(u) $ is written in the basis elements 
$ T_w E_{A} $ of $G$, see (\ref{basis}). Here by abuse of notation, we write $ {\bf n}$ for the 
unique maximal set partition in $ {\cal P}_n $. Its only block is $ {\bf n} $.

\medskip
\noindent
With this we may construct a 
bilinear form $ \langle \cdot, \cdot \rangle $ on ${\cal E}_n(u)$ by
$$ \langle x , y \rangle  = \epsilon ( x^{\ast} y )  \,\,\,\,\,\, \mbox{ for } \, x, y \in {\cal E}_n(u) $$
And then we finally obtain a homomorphism $ \varphi$ by the rule 
$$ \varphi:  {\cal E}_n(u) \rightarrow  {\cal E}_n(u)^{\prime}: x \mapsto 
( y \mapsto \langle x, y \rangle )$$

\begin{Thm} With the above definitions, we get that
$ \varphi $ is an isomorphism of left
${\cal E}_n(u) $-modules.
\end{Thm} 
\begin{pf*}{Proof}
One first checks that the bilinear form satisfies 
$$ \langle x y, z \rangle = \langle  y, x^{\ast} z \rangle  \,\,\,\,\,\, \mbox{ for all } x,y,z \in {\cal E}_n(u) $$
which amounts to saying that $ \varphi $ is ${\cal E}_n(u) $-linear.

\medskip Since ${\cal E}_n(u) $ is finite dimensional, it is now
enough to show that $ \langle \cdot , \cdot \rangle $ is
nondegenerate. For this we first observe that our construction of $
\langle \cdot , \cdot \rangle $ is valid over $ \A $ as well and hence
also defines a bilinear form $ \langle \cdot , \cdot \rangle_{\A} $ on
$ {\cal E }^{\A}_n(u) $. It is enough to show that $ \langle \cdot ,
\cdot \rangle_{\A} $ is nondegenerate. 
Suppose $ a \in {\cal E }^{\A}_n(u) $. 
Then as in the proof of theorem {\ref{Gbasis}} we can write it in the form
$ a = (u-1)^N a^{\prime}$ where $ a^{\prime} = \sum_{ g \in G }
\, \lambda_g \, g $ and where $ \lambda_g(1) \not= 0$ for at least one $ g $.
Letting $ \pi:  {\cal E}^{\A}_n(u)  \rightarrow  {\cal E}_n(1) $ be the specialization 
map we
have $ \pi(a^{\prime}) \not= 0 $ since it was shown in the proof of that theorem that 
$ G $ is a basis of $ {\cal E}^{\A}_n(1) $ as well.

\medskip
Let us denote by $ \langle \cdot , \cdot \rangle_1 $ the bilinear form on $ {\cal E}_n(1) $ constructed similarly
to $ \langle \cdot , \cdot \rangle $. Then we 
have that 
$$ \langle \pi(a), \pi(b) \rangle_1 = 
\langle a , b \rangle_{\A} \otimes_{\A} \C  \,\,\,\,\,\, \mbox{   for all } a,b \in  {\cal E}^{\A}_n(u) $$
since $ \pi $ is multiplicative and satisfies $ \pi( a^{\ast}) = \pi(a)^{\ast} $.
We are now reduced to proving that $ \langle \cdot , \cdot \rangle_1 $ is nondegenerate. 
Let us therefore consider an arbitrary $ a = \sum_{w,A } \lambda_{w, A }  E_A T_w  \in {\cal E}_n(1) $, where $ \lambda_{w, A } \in \C $.
Let $ A_0 \in {\cal P}_n $ be minimal 
subject to the condition that $ \lambda_{w,A_0} \not= 0 $ for some $ w $.
Take $ z  \in S_n $ with $ \lambda_{z, A_0} \not= 0 $ and define 
$$ b = E_{A_0} \prod_{ A_0 \subsetneq A} (1-E_{A}) \, T_{z} $$
We claim that $ \langle b, a \rangle_1 \not= 0 $. Indeed, since $ u= 1 $ we have 
$$ b^{\ast}a = T_{z}^{-1}  \prod_{ A_0 \subsetneq A} (1-E_{A}) E_{A_0} a $$
Since $ A_0 $ was chosen minimal, there can be no cancellation of the coefficient of $E_{A_0} T_{z} $ in $ E_{A_0} a $ 
which hence is $ \lambda_{z, A_0}  $. All $ E_A $ appearing in the expansion of 
$ E_{A_0} a $ with respect to 
the basis $ E_A T_w $ satisfy $ A_0 \subseteq A $. Except for $ E_{A_0} $ they are 
all killed by $ \prod_{ A_0 \subsetneq A} (1-E_{A}) $. By this we get 
$$ T_{z}^{-1} \prod_{ A_0 \subsetneq A} (1-E_{A}) E_{A_0} a =  \lambda_{z, A_0}  \,
 T_{z}^{-1} \prod_{ A_0 \subsetneq A} (1-E_{A}) E_{A_0}  T_{z} $$
The coefficient of $ E_{\bf n }$ in this expression is by corollary \ref{setreflection} 
equal to the 
coefficient of $ E_{\bf n }$ in 
$$ \lambda_{z, A_0}  \prod_{ A_0 \subsetneq A} (1-E_{A}) E_{A_0}   $$
On the other hand, the coefficient of $ E_{\bf n }$ in  $ \prod_{ A_0 \subsetneq A} (1-E_{A}) E_{A_0}   $
is given by the Moebius function associated with the partial order $ \subset $ on $ {\cal P}_n $. 
It is equal to $ (-1)^{k-1} k! $, where $k $ is the number of 
blocks of $ A_0 $. Summing up we find that $ \langle b, a \rangle_1 = (-1)^{k-1} \lambda_{z, A_0}  k! \not= 0  $
which proves the theorem.
\end{pf*}

\section{Classification of the irreducible representations}
In this section we give the classification of the irreducible representations of $ {\cal E}_n(u) $.

\medskip
For $ M $ a left $ {\cal E}_n(u) $-module we make its linear dual $ M^{\prime} $ into a left $ {\cal E}_n(u) $-module 
using the antiautomorphism $ \ast $. If $ M $ is a simple $ {\cal E}_n(u) $-module 
then any $ m \in M \setminus \{0\}$ 
defines a surjection 
$$ {\cal E}_n(u) \rightarrow M, \, x \mapsto x m \,\,\,\,\,\,\,\,\,  \mbox{    for} \,x \in {\cal E}_n(u) $$
By duality and by the last section, we then get an injection of $ M^{\prime} $ into ${\cal E}_n(u) $. On the other hand, the canonical 
isomorphism $ M \rightarrow M^{\prime \prime} $ is $ {\cal E}_n(u) $-linear because 
$ \ast \ast = Id $ and so we conclude that 
all simple $ {\cal E}_n(u) $-modules appear as left ideals in $ {\cal E}_n(u) $. 

\medskip
Let now $ I $ be a simple left ideal of $ {\cal E}_n(u) $ and let $ x_0 \in I \setminus \{0 \}  $. 
Since the tensor space $ V^{ \otimes n }$ is a faithful $ {\cal E}_n(u) $-module, we find a $ v \in 
 V^{ \otimes n }$ such that $ x_0 v \not= 0 $. But then the $ {\cal E}_n(u) $-linear map
$$ I \rightarrow V^{ \otimes n }, \, x \mapsto x v \,\,\,\,\,\,\, \mbox{  for } x \in I $$ 
is nonzero, and therefore injective. We conclude that all simple $ {\cal E}_n(u) $-modules appear as submodules 
of $ V^{ \otimes n }$.

\medskip
Consider a simple submodule $ M $ of $ V^{ \otimes n } $. Take $ A_0 \subset { \bf n}   $ 
maximal such that $ E_{A_0} M \not= 0$. 
By section 3, in the two extreme situations $ A_0 = \emptyset $ or $ A_0 = {\bf n} $ we can give a precise description of 
$ M $, since in those cases $ M$ is a module for $   K S_n $ 
or $ H_n(u) $. 
In other words, $ M $ is the pullback of a Specht module $ S(\lambda )$ for $ K S_n $ or
a Specht module $ S_u(\lambda )$ for $ H_n(u) $ as
described in section 3. The general case is going to be a mixture of these two cases as we shall 
explain in this section.

\medskip
Let $ {\cal L}_n  $ be the set of tuples 
$$  {\cal  L }_n = \{ \, ( \lambda^{s} , m_s  ,  \mu^s ) \, | \, s = 1, \ldots , k \,  \} $$
where $ \lambda^{s} $ is a partition, $ m_s $ a positive integer and $ \mu^s $ a 
partition of $ m_s $ such that 
$ \sum_s  m_s \, | \lambda^{{s}} | = n  $ and such that $ \lambda^1 < \lambda^2 <
\ldots < \lambda^k$ where $ < $ is the 
total order on partitions defined above.

\medskip
Suppose $ \Lambda = ( \lambda^{s} , m_s  ,  \mu^s ) \in {\cal L}_n $. We associate to it the vector 
$ v_{\Lambda} \in V^{\otimes n } $ defined in the following way 
$$ v_{\Lambda} := v_{\lambda^{1}}^1 \otimes v_{ \lambda^{1}}^2 \otimes \ldots \otimes 
v_{\lambda^{2}}^{m_1 +1  } \otimes v_{\lambda^{2}}^{m_1+ 2} \otimes \ldots  \otimes
v_{\lambda^{k}}^{l}
$$
where $ l := \sum_s m_s $ and where 
for any integer partition (even composition) $ \mu 
= (\mu_1, \mu_2, \ldots, \mu_r ) $ of $  m $ and any integer $ i $ we define 
$ v_{\mu}^i  \in V^{\otimes m } $ as follows 
$$ v_{\mu}^i := (v_1^i)^{ \otimes  \mu_1} \otimes (v_2^i)^{\otimes \mu_2} \otimes 
\ldots \otimes (v_r^i)^{\otimes \mu_r}  $$

\medskip
We moreover associate to $  \Lambda = (\lambda^{s} , m_s  ,  \mu^s ) $ 
the set partitition $ A_{ \Lambda} \in {\cal P}_n $, 
that has blocks of consecutive numbers, 
the first $ m_1 $ blocks being of size $ | \lambda^1|  $, the next $ m_2 $ blocks of size 
$ | \lambda^2 |$ and so on. The blocks correspond to the factors of $ v_{\Lambda} $ that 
have equal upper indices.
Note that it is possible that $ | \lambda^1| = | \lambda^2| $ 
making the first $ m_1 + m_2 $ blocks of equal size and so on. 
Writing $ A_{ \Lambda} = (I_1, I_2, \ldots, I_l ) $ we set 
$$ \begin{array}{c}
S_{\Lambda} :=  S_{m_1} \times S_{m_2} \times  \ldots \times S_{m_k } \\
H_{\Lambda}(u) := H_{I_1}(u) \otimes H_{I_2}(u) \otimes \ldots  \otimes H_{I_l}(u)
\end{array}
$$
Let $ \iota_j $ be the group isomorphism from 
$ S_{m_j } $ to $ 1 \times \ldots \times S_{ m_j } \times \ldots \times 1 $ and 
also the 
algebra isomorphism from $ H_{I_j}(u) $ to 
$ 1 \otimes \ldots \otimes H_{I_j}(u)  \otimes \ldots  \otimes 1 $. 

Corresponding to $ A_{ \Lambda} $ there is an analogous block
decomposition of the factors of $ V^{\otimes n } $ and $ S_{\Lambda} $ acts on this 
by permutation of the blocks. 

\medskip
Let us illustrate this action on an example.
Take $ n = 6$, $ k = 1 $ and $ \Lambda = ( \lambda, 2, \mu) $ where 
$ \lambda= (2,1) $ and $ \mu = (1,1) $. 
Then $ A_{\Lambda} = \{ (1,2,3), (4,5,6) \} $ and $ S_{ \Lambda} $ is the group of order two
that permutes the two blocks, thus generated by $ \sigma 
= (1,4) (2,5) (3,6) $.
In other words
$$ v_{\Lambda} = v_1^1 \otimes v_1^1 \otimes v_2^1 \otimes v_1^2 \otimes v_1^2 \otimes v_2^2 
\mbox{  and  } 
\sigma v_{\Lambda} = v_1^2 \otimes v_1^2 \otimes v_2^2 \otimes v_1^1 \otimes v_1^1 \otimes v_2^1 
$$
In general, we have that 
\begin{equation} T_{\sigma} v_{\Lambda} = \sigma v_{\Lambda} \mbox{    for } 
 \sigma \in S_{\Lambda} 
\end{equation}
since 
in a reduced expression $ \sigma = \sigma_{i_1} \sigma_{i_2} \ldots \sigma_{i_N} $ 
the action of each $ \sigma_{i_j} $ and $ T_{i_j} $ on $ v_{\Lambda} $ will only involve 
distinct upper indices.

\medskip
In the above example, we have  
$ \sigma = \sigma_3 \sigma_4 \sigma_5 \sigma_2 \sigma_3 \sigma_4  \sigma_1 \sigma_2 \sigma_3  
\in S_{\Lambda} $
and hence 
$$ T_{\sigma} = T_3 T_4 T_5 T_2 T_3 T_4  T_1 T_2 T_3 \in  {\cal E}_n(u) $$
Both $ \sigma $ and $ T_{\sigma}$ will move the first $ v_1^2 $ to the first position, then
the second $ v_1^2 $ to the second position and finally $ v_2^2 $ to the third position.

\medskip
We consider the row and column (anti)symmetrizer $ r_{\mu^i},  c_{\mu^i} \in K S_{ | \mu^i  | }$
of the partitions $ \mu^i$ as elements of 
$ {\cal E}_n(u) $ by mapping each occurring $ \sigma $ to $ T_{\iota_i(\sigma)} $.
By corollary \ref{setreflection}, we then get that 
$ r_{\mu^i} $ and $ c_{\mu^i} $ commute with $ E_{A_{\Lambda}} $.

\medskip
We define 
$ w_{\Lambda}   := ( r_{\mu^1} \otimes r_{\mu^2} \otimes \ldots \otimes r_{\mu^k}  ) v_{\Lambda} $.
It has the form 
$ w_{\Lambda} := w_{\lambda_1}^{ \mu_1 } \otimes \ldots \otimes w_{\lambda_k}^{ \mu_k } $
where we for general $ \lambda, \mu $ define 
$$ w_{\lambda}^{ \mu } := \sum_{ \sigma \in r_{\mu}} v_{ \lambda}^{ \sigma(1)} \otimes \ldots 
\otimes  v_{ \lambda}^{ \sigma(m)}         $$ 
where $ | \mu | = m $.
We define the 'permutation module' as $$ M(\Lambda) := {\cal E}_n(u)
= {\cal E}_n(u) w_{\Lambda} $$ 
Define now 
$$ e_{\Lambda } := 
( c_{\mu^1} \otimes c_{\mu^2} \otimes \ldots \otimes c_{\mu^k}  )
(c_{\lambda^1}(u)^{\otimes m_1} \otimes \ldots \otimes c_{\lambda^k}(u)^{ \otimes m_k} )
E_{A_{\Lambda}}$$
where $ c_{\lambda^i}(u) $ is as in section 4.
Note that the three factors of $ e_{ \Lambda} $ commute by the definitions and 
corollary \ref{setreflection}. 
We define the 'Specht module' as 
$$ S(\Lambda) := {\cal E}_n(u)e_{\Lambda}  w_{\Lambda} \subset M(\Lambda) $$
Actually, the factor $ E_{A_{\Lambda} } $ could have been left out of $ e_{\Lambda} $ 
in the definition of the Specht module, 
since it commutes with $  r_{\mu^1} \otimes r_{\mu^2} \otimes \ldots \otimes r_{\mu^k}   $ 
and $ E_{A_{\Lambda} } w_{ \Lambda} = w_{ \Lambda} $ by the next lemma \ref{E-action}, 
but for later use we prefer
to include it in $ e_{\Lambda} $.

\begin{lemma} In the above setting we have that 
\begin{equation}{\label{E-action}}
E_B w_{\Lambda} = \left\{ \begin{array}{cc}
w_{\Lambda} & \mbox{if } B \subseteq A_{\Lambda} \\
0  & \mbox{otherwise }
\end{array}
\right.
\end{equation}
\end{lemma}
\begin{pf*}{Proof} 
If $ B \subseteq A_{\Lambda} $ 
this is an immediate consequence of the definitions. 
If $ B \not\subseteq A_{\Lambda} $ there are $ i,j \in {\bf n} $ belonging to the 
same block of $ B $ and to different blocks of $ A_{\Lambda} $, let these be 
$ I_{\alpha(i)} $ and $  I_{\alpha(j)} $.  
Since $ E_{ij} $ is a factor of $ E_B$ it is enough 
to show that $ E_{ij} \sigma v_{\Lambda} = 0 $ for $ \sigma \in S_{\Lambda} $.
But from formula (\ref{innermost}) we have that 
$$ E_{ij} = T_{j-1} T_{j-2} \ldots T_{i+1}  E_{i}  
T_{i+1}^{-1} \ldots T_{j-2}^{-1} T_{j-1}^{-1} $$
Using it we can decompose $ E_{ij} $ from the right to the left 
in an element of $ \iota_{\alpha(j)} (H_{I_{\alpha(j)}}) $, 
followed by the product of the remaining $ T_{k}^{-1} $, then $ E_{i} $ 
and finally the product of the $ T_k $. 
The action of $ \iota_j (H_{I_{\alpha(j)}}) $ on $ \sigma v_{\Lambda} $ produces a linear 
combination of basis elements $ v $ of $ V^{\otimes n } $ where all appearing $ v $ 
are obtained from $ \sigma v_{\Lambda} $ by permuting the factors corresponding to 
the block $ I_{\alpha(i)} $. The upper indices of the factors of $ v $ are exactly as 
those of $\sigma  v_{\Lambda}$.
The product of $ T_{k}^{-1} $ 
acts on each $ v $ by permuting the first factor of the $ I_{\alpha(j)} $
block to the $i+1$st position, that is inside the $ I_{\alpha(i)} $ block. 
But $ E_i $ acts as zero on this and the lemma follows.
\end{pf*}

The main result of this section is the following theorem.
\begin{Thm}{\label{simplesubmodule}} $S(\Lambda)$ is a simple module for $ {\cal E}_n(u) $. 
The simple ${\cal E}_n(u) $-modules are classified by $ S(\Lambda) $ for $ \Lambda \in {\cal L}_n$.
\end{Thm} 
\begin{pf*}{Proof} 
Write for simplicity $ A:= A_{\Lambda} $. 

\medskip
Our first step is to show that $ e_{\Lambda} M(\Lambda ) = K  e_{\Lambda} w_{\Lambda} $.
For this we take $ x \in {\cal E}_n(u) $ and first consider the element 
$ E_{A} x w_{ \Lambda} \in M(\Lambda) $.

\medskip
We can write $ x$ as a linear combination of elements $ E_B T_w $ from our basis $ G $.
By 
corollary {\ref{inclusion}}, $ E_A E_B $ is 
equal to a $ E_C $ for $ C $ with $ A \subseteq  C $. By lemma \ref{E-action} and 
corollary \ref{setreflection} we have that 
$ E_C  T_w w_{ \Lambda}  = T_w  E_{ w^{-1}C }  w_{ \Lambda}  = 0 $ unless 
$ w^{-1}C = A $, since $ A \subseteq  C $. We may therefore 
assume that $ B = A $ and $ A = w A $ such that 
$ E_A x $ is a linear combination of elements of the form 
$ E_A T_w $ where 
$ T_w $ permutes the blocks of $ A $ of equal cardinality.

\medskip
Thus, let $ \overline{S_{{\Lambda}}} \le S_n $ be the subgroup consisting of the permutations of 
the blocks of $ A $ of equal cardinality. Note that 
$ S_{ \Lambda} \le  \overline{S_{\Lambda}} $, the inclusion being strict in general. 
As in the case of $ S_{ \Lambda} $, the 
elements of $  \overline{S_{\Lambda}} $ can be seen as elements of  
$ {\cal E}_n(u) $, by the map $ z \mapsto T_z $.

\medskip
In this notation, if $ E_A x w_{\Lambda} $ is nonzero it is 
a linear combination of elements of the form
\begin{equation}{\label{summingup}}
T_z \, ( T_{w_{1}}  \otimes   T_{w_{2}} \, 
 \otimes \ldots \otimes
T_{w_{l}}  ) w_{\Lambda}
\end{equation}
where $ z \in \overline{S_{\Lambda}}$ and 
$ T_{w_{1}} \otimes   T_{w_{2}} \, 
 \otimes \ldots \otimes
T_{w_{l}}  \in  H_{ \Lambda }(u) $ and where we used that $ E_A $ commutes with the 
other factors and $ E_{A} w_{ \Lambda} = w_{ \Lambda} $.
Since the upper indices of the $ w_{ \lambda^i}^j $ are distinct, $ T_z $ acts by 
permuting the $ T_{w_i}$-factors.

\medskip
We need to show that $ z \in S_{\Lambda}$ and therefore consider the action on 
$  c_{\lambda^1}(u)^{\otimes m_1} \otimes \ldots \otimes c_{\lambda^k}(u)^{ \otimes m_k}  $ 
on (\ref{summingup}). 
Let from this $ \lambda^1, \lambda^2, \ldots, \lambda^t $ be the partitions 
with $ | \lambda^i | =  | \lambda^1 | = | I_1 | $. Note that in general $ t \geq m_1 $. 
Since the $ \lambda^i $ are ordered increasingly, we get by lemma \ref{futurereference} that the product is nonzero only if each 
factor
$ c_{\lambda^k}(u) $ 
 of $ c_{\lambda^1}(u)^{ \otimes m_1}  \otimes \ldots \otimes c_{\lambda^t}(u)^{ \otimes m_t}   $
acts in a 
$ T_{w_a}\,  v_{\lambda^{k}}^{\, \sigma(a)} $-factor of (\ref{summingup}), i.e. a factor with the same $ \lambda^k $ appearing as index.
This argument extends to the other factors of 
$  c_{\lambda^1}(u)^{\otimes m_1} \otimes \ldots \otimes c_{\lambda^k}(u)^{ \otimes m_k}  $ 
and so we may assume that $ z \in S_{\Lambda}$ as claimed.


\medskip
After this preparation, we can show the claim about $ e_{\Lambda} M(\Lambda ) $. 
We take $ x \in {\cal E}_n(u) $ and consider $ e_{\Lambda} x w_{ \Lambda} $. 
By the above, 
it is a linear combination of elements of the form 
$$   
( c_{\mu^1} \otimes \ldots \otimes c_{\mu^k}  )T_z(c_{\lambda^1}(u) \otimes     \ldots \otimes 
c_{\lambda^l}(u) )  
( T_{w_{1}}  	
 \otimes \ldots \otimes
T_{w_{l}}) 
w_{ \Lambda}
$$
where 
$ T_{w_{1}} \otimes   T_{w_{2}} \, 
 \otimes \ldots \otimes
T_{w_{l}}  \in  H_{ \Lambda }(u) $
and where $ z \in S_{\Lambda}$ such that $ T_z $ commutes with 
$ c_{\lambda^1}(u) \otimes     \ldots \otimes 
c_{\lambda^l}(u)$. 
We now use the formulas (\ref{schur_hecke}), (\ref{schur_sym}) and the definition of 
$ w_{\Lambda} $ 
to rewrite this as 
$$ 
\begin{array}{r}
C_1 ( c_{\mu^1} \otimes \ldots \otimes c_{\mu^k}  )T_z(
s_{\lambda^1}(u) \otimes     \ldots \otimes 
s_{\lambda^l}(u) )  
w_{ \Lambda} = \\
C_2 ( c_{\mu^1} \otimes \ldots \otimes c_{\mu^k}  )T_z(
s_{\lambda^1}(u) \otimes     \ldots \otimes 
s_{\lambda^l}(u) )  
( r_{\mu^1} \otimes \ldots \otimes r_{\mu^k}  )w_{ \Lambda} = \\
C_2 ( c_{\mu^1} \otimes \ldots \otimes c_{\mu^k}  )T_z
( r_{\mu^1} \otimes \ldots \otimes r_{\mu^k}  )
(s_{\lambda^1}(u) \otimes     \ldots \otimes 
s_{\lambda^l}(u) )  
w_{ \Lambda} = \\
C_3 ( s_{\mu^1} \otimes \ldots \otimes s_{\mu^k}  )
(s_{\lambda^1}(u) \otimes     \ldots \otimes 
s_{\lambda^l}(u) )  
w_{ \Lambda} =  \\
C_4 ( c_{\mu^1} \otimes \ldots \otimes c_{\mu^k}  ) 
(c_{\lambda^1}(u) \otimes     \ldots \otimes 
c_{\lambda^l}(u) )  
w_{ \Lambda}  = C_4 e_{ \Lambda} w_{ \Lambda} 
\end{array}
$$
where the $ C_i \in K$ are constants and 
where 
we used that 
$  r_{\mu^1} \otimes \ldots \otimes r_{\mu^k}   $ 
 commutes with  
$ c_{\lambda^1}(u) \otimes     \ldots \otimes 
c_{\lambda^l}(u)   $ and $ r_{\lambda^1}(u) \otimes     \ldots \otimes 
r_{\lambda^l}(u)   $ since $ r_{\mu^1} $ permutes over equal factors $ c_{\lambda^1}(u) $ etc.
For $ z =1 $ all the constants are nonzero since the Young symmetrizers $ s_{\lambda}(u) $ and 
$ s_{\mu} $ are idempotents up to nonzero scalars and we have then finally proved 
that $ e_{\Lambda} M(\Lambda ) = K e_{\Lambda } w_{\Lambda } $, as claimed. Since $ S(\Lambda ) \subseteq  M(\Lambda)  $ we also have  
$ e_{\Lambda} S(\Lambda ) \subseteq  K e_{\Lambda } w_{\Lambda } $.

\medskip
We now proceed to prove that $ S(\Lambda ) $ is a simple module for $  {\cal E}_n(u) $. 
We do it by setting up of version of James's submodule theorem, [Ja].
Assume therefore $ N \subset S(\Lambda)  $ is a submodule. If $ e_{\Lambda} N \not= 0 $, we 
have by the above that $ e_{\Lambda} N $ is a scalar multiple of 
$e_{\Lambda} w_{\Lambda} $ and so $ N = S(\Lambda ) $.

\medskip
In order to treat the other case $ e_{\Lambda} N = 0 $, we 
define a bilinear form on $ V^{\otimes n } $ 
by setting 
$$ \langle v_{i_1 }^{j_1} \otimes \ldots \otimes v_{i_n }^{j_n}, v_{i_1^{\prime} }^{j_1^{\prime}} \otimes \ldots \otimes 
v_{i_n^{\prime} }^{j_n^{\prime}} \rangle = v^{\overline{i}} \delta_{\overline{i}= 
  \overline{i^{\prime}}, \overline{j}=\overline{j^{\prime}}}$$
and extending linearly, where 
we write $ \overline{i} = (i_1, i_2 , \ldots , i_n) $ and similarly for
$ \overline{i^{\prime}}, \overline{j},\overline{j^{\prime}}$.
The power 
$ v^{\overline{i}} $ is defined as follows. 
Order $ v_{i_1 }^{j_1} \otimes \ldots \otimes v_{i_n }^{j_n} $ 
by first moving all factors $ v_{i_k }^{j_k} $ with minimal upper indices 
to the left of $ v_{i_1 }^{j_1} \otimes \ldots \otimes v_{i_n }^{j_n} $
but maintaining their relative position, then moving the factors $ v_{i_k }^{j_k} $ 
with second smallest upper indices to the positions just to the right of the first ones and so on.
This gives a permutation $ \sigma \in S_n $ such that $ \sigma (
v_{i_1 }^{j_1} \otimes \ldots \otimes v_{i_n }^{j_n}) $ has increasing upper indices, let these be
$ f(1), f(2) , \ldots , f(m) $ without repetitions. 
We then find compositions $ \tau_i, i = 1, \ldots, m $ and minimal coset representations
$ w_i \in S_{| I_{\tau_i}|}/S_{ \tau_i} $ 
such that 
$$ \sigma
(v_{i_1 }^{j_1} \otimes \ldots \otimes v_{i_n }^{j_n}) = 
 {w_{1}}\,  v_{\tau^{1}}^{ f(1)} \otimes   {w_{2}} \, 
v_{\tau^{2}}^{f(2)} \otimes \ldots \otimes
{w_{m}} \, v_{\tau^{m}}^{f (m)} 
$$
and define $ v^{\overline{i}} := v^{ \sum l(w_i)} $. 

This bilinear form is modelled on the one for the tensor space module 
for Hecke algebras, [DJ], and inherits 
from it the following invariance property
$$ \langle x v , w\rangle = \langle v , x^{\ast} w \rangle \, \, \, \mbox{  for all} \, \,  x \in  {\cal E}_n(u), v, w \in V^{\otimes n }$$
where $ \ast $ is as in section 4. 
We have that 
$$
c_{\lambda}^{\ast} = c_{\lambda},  \, \, \,   r_{\lambda}^{\ast} = r_{\lambda},  \,\,\,  
c_{\lambda}(u)^{\ast} = c_{\lambda}(u),  \,\,\, 
 r_{\lambda}(u)^{\ast} = r_{\lambda}(u)
$$ 
where we used that $  \ast $ is an antiautomorphism to show for instance that
$ T_{ w_{\lambda}^{-1} } y_{\lambda^{\prime} } T_{ w_{\lambda} }^{\ast} = 
T_{ w_{\lambda}^{-1} } y_{\lambda^{\prime} } T_{ w_{\lambda} } $.
Since the factors of 
$ e_{\Lambda} $ commute, we also have that
$$ e_{\Lambda}^{\ast} = e_{\Lambda}  $$
We are now in position to finish the treatment of the case $ e_{\Lambda} N = 0 $. 
We have
$$ 0= \langle e_{\Lambda} N, M(\Lambda) \rangle = \langle  N, e_{\Lambda} M(\Lambda) \rangle  = 
\langle  N, e_{\Lambda} w_{\Lambda}  \rangle  $$ 
which implies that $ \langle  N, S(\Lambda )  \rangle = 0 $
that is $ N \subset S(\Lambda )^{\perp} $. 
Since $ u $ is generic, we have that $ \langle e_{\Lambda} w_{\Lambda},e_{\Lambda} w_{\Lambda} \rangle
\not=0 $ and therefore 
$ S(\Lambda )  \cap S(\Lambda )^{\perp} =  0 $. This gives a contradiction unless $ N= 0 $. 
We have therefore proved that $ S(\Lambda) $ is simple.

\medskip 
We next prove that different choices of parameters give 
different modules
$ S(\Lambda)  $. Take $ \Lambda $ as before 
and suppose $ { \Upsilon } = 
(  ({ \nu}^t ),\,( {n}_t) , \, ({\tau}^t))  \in {\cal L}_n $ 
such that $ S(\Lambda) \cong S({\Upsilon}) $.
The element $ A \in {\cal P}_n $ associated with $ S(\Lambda) $ is maximal with respect to having blocks of consecutive numbers such that
$ E_A S(\Lambda) \not=  0 $. Hence, if $ B \in {\cal P}_n $ is the 
element associated with $ S({\Upsilon}) $, we have that 
$ A = B$. But then $ ({ \lambda}^s ) $ and $ ({ \nu}^t ) $ must be partitions of the same numbers, 
corresponding to the block sizes of 
$ A $, or $ B$. Both $ c_{\lambda^1}(u) \otimes     \ldots \otimes 
c_{\lambda^l}(u)   $ and 
$ c_{\nu^1}(u) \otimes     \ldots \otimes 
c_{\nu^l}(u)   $ act nontrivially in $ E_A S(\Lambda) $ 
and hence by lemma \ref{futurereference} we have $ \lambda^i \leq \nu^i $ and 
$ \lambda^i \geq \nu^i $ that is $ \lambda^i = \nu^i $. 
Similarly, we get $ ({ \mu}^s )=  (  \tau^t ) $. This proves the claim.

\medskip It remains to be shown that any simple module $ L $ is of the form 
$ S(\Lambda ) $ for some $ \Lambda \in {\cal L}_n $.
We saw in the remarks preceding the theorem, that 
it can be assumed that $ L \subset  V^{\otimes n } $. Choose 
$ A = \{ I_1, \ldots, I_l \} 	\in {\cal P}_n $ maximal 
with respect to having blocks of consecutive numbers and 
$ E_A L \not=  0 $. For $ \sigma \in S_n $, the map 
$ \varphi^{\sigma}: V^{\otimes } \rightarrow V^{\otimes } $ 
given by 
$$ \varphi^{\sigma} : v_{i_1}^{j_1} \otimes \ldots \otimes v_{i_n}^{j_n} \rightarrow 
v_{i_1}^{\sigma(j_1)} \otimes \ldots \otimes v_{i_n}^{\sigma(j_n)} $$
is an $ {\cal E}_n(u) $-linear isomorphism and replacing $ L $ by $ \varphi_{\sigma} L $ for 
an appropriately chosen $ \sigma $ we may assume that 
$ | I_i |\leq | I_{i+1} |$ for all $ i $.
We have now that $ E_A L $ is a module for the tensor product 
$ H_{I_1}(u)  \otimes \ldots \otimes H_{I_l}(u) $.
Choose for each $ I_i $ a partition $ \lambda_i $ of $ | I_i | $ 
such that the 
product $ c_{\lambda^1}(u) \otimes c_{\lambda^2}(u) \otimes \ldots \otimes c_{\lambda^l}(u) $ 
acts nontrivially in $ E_A L $. Choose next partitions $ \mu^i $ such that 
$ s_{\mu^1} \otimes s_{\mu^2} \otimes \ldots  s_{\mu^k} $ acts nontrivially in 
$  (c_{\lambda^1}(u)  \otimes  \ldots \otimes c_{\lambda^l}(u)) E_A L $.  
The data so collected give rise to a $ \Lambda $ 
with $ S(\Lambda ) = {\cal E}_n(u) e_{\Lambda} w_{\Lambda} \subset L $. 
But since $ L $ is simple, the inclusion must be an equality.
With this we have finally proved all statements of the theorem.
\end{pf*}

\medskip

\medskip
Let us work out some low dimensional cases. 
For $ n =2 $ we have the following possibilities for $\Lambda$:
$$\begin{array}{c} 
(\lambda^1,m_1, \mu^1)=(\,  {\small{\tableau[scY]{, |
}}}\,  , 1,\,  {\small{\tableau[scY]{ |
}}} \,  ), \, \,\,\,
(\lambda^1,m_1, \mu^1)=(\,  {\small{\tableau[scY]{ | |
}}}\,  , 1,\,  {\small{\tableau[scY]{ |
}}} \,)     \\
(\lambda^1,m_1, \mu^1)=(\,  {\small{\tableau[scY]{ | 
}}}\,  , 2,\,  {\small{\tableau[scY]{ ,
}}} \,  ), \, \,\,\,
(\lambda^1,m_1, \mu^1)=(\,  {\small{\tableau[scY]{  |
}}}\,  , 2,\,  {\small{\tableau[scY]{ | | 
}}} \,)  

\end{array}
$$
They all give rise to irreducible representations of dimension one. 
The first two are the one dimensional representations of $H_2(u) $. By our construction  
the third is given by $ v_1^1 \otimes v_1^2 + v_1^2 \otimes v_1^1 $ and the last by 
$ v_1^1 \otimes v_1^2 - v_1^2 \otimes v_1^1 $. They correspond to the trivial and the 
sign representation of $ K S_2 $. The square sum of the 
dimensions is 4, which is also the dimension of $ {\cal E}_2(u) $.

\medskip
For $n=3 $ we first write down the multiplicity free possibilities of $\Lambda $, i.e. those having $ m_s = 1 $ and so $ \mu_s 
= {\small{\tableau[scY]{ |}}} \, $ for all $ s $. They are 
$$ 
\begin{array}{lll} 
(\lambda^1)  = ( \, {\small{\tableau[scY]{,, |}}} \, ) , & (\lambda^1 )  =  ( \, {\small{\tableau[scY]{ , | |}}} \, ),    &
(\lambda^1 ) = ( \, {\small{\tableau[scY]{| | |  }}} \, )  \\  ( \, \lambda^1, \lambda^2 \, ) =  ( {\, \small{\tableau[scY]{ |}}  , \small{\tableau[scY]{, |}}\, }\,  ),  & 
(\, \lambda^1, \lambda^2 \, ) =(\,  {\small{\tableau[scY]{   |}}} , \,  {\small{\tableau[scY]{|  |}}} \,)   &
\end{array}
$$
The first three of these are the Specht modules for $ H_3 (u)  $, their dimensions are respectively 1,2 and 1. 
The fourth is given by the vector
$  v_1^1 \otimes v_1^2  \otimes  v_1^2 $ and the last by the vector 
$  v_1^1 \otimes (v_1^2 \otimes v_2^2 - u^{-1} v_2^2\otimes v_1^2)  $, according to our construction.
In both cases, one gets dimension three.

\medskip
Allowing multiplicities, we have the following possibilities:
$$\begin{array}{cc} 
(\lambda^1, m_1, \mu^1) = ( \,  {\small{\tableau[scY]{  |
}}}\,, 3 , \,  {\small{\tableau[scY]{  | | |
}}}\,), &
(\lambda^1, m_1, \mu^1) = ( \,  {\small{\tableau[scY]{  |
}}}\,, 3 , \,  {\small{\tableau[scY]{ , | |
}}}\,) 
\end{array} 
$$
and 
$\begin{array}{c} 
(\lambda^1, m_1, \mu^1) = ( \,  {\small{\tableau[scY]{  |
}}}\,, 3 , \,  {\small{\tableau[scY]{  | | | 
}}}\,) 
\end{array}
$.
We get the Specht modules of $ K S_3 $ of dimensions 1,2 and 1. 

\medskip
The square sum of all the dimensions is 30, in accordance with the dimension of $ {\cal E}_3(u) $. We have thus proved that 
$ {\cal E}_n(u) $ is semisimple for $ n = 2 $ and $ n= 3 $. 

\medskip
The classification of the simple modules for $ n= 2 $ and $ n= 3 $ has also been done in [AJ] with a different method.

\section{Questions}
The results of the paper raise a number of questions. 

\medskip
There is a canonical inclusion
$ {\cal E}_n(u) \subset  {\cal E}_{n+1}(u) $ which at diagram level is given by adding a
through line to the right of a diagram element from $ {\cal E}_n(u) $. 
It gives rise to restriction 
and induction functors $ \res $ and $ \ind $, that should obey a  
branching rule 
for 
the decomposition of $ \res S(\Lambda) $. Our first question is to 
give a description of it.
Apart from the independent interest in such 
a branching rule, one possible application would 
be to obtain a dimension formula for $ S(\Lambda) $.

\medskip
We do not know what the general branching rule looks like, 
but using the above calculations, we can at least explain the cases $ n=2,3 $, 
corresponding to 
$ \res S(\Lambda) $ for $ \Lambda \in {\cal P}_{2} $ and 
$ \Lambda \in {\cal P}_{3} $,
These cases are rather easy, since one only needs consider $ n = 3 $, $ \Lambda = ( \lambda^s, m_s, \mu^s) $, 
$ m_s = 1 $ and $ \mu^s $ trivial and 
$$ \begin{array}{lll} 
    ( \, \lambda^1, \lambda^2 \, ) =  ( {\, \small{\tableau[scY]{ |}}  , \small{\tableau[scY]{, |}}\, }\,  ),  & 
(\, \lambda^1, \lambda^2 \, ) =(\,  {\small{\tableau[scY]{   |}}} , \,  {\small{\tableau[scY]{|  |}}} \,)   &
\end{array}
$$
because, as we saw above, all 
other choices of $ \Lambda $ give Specht modules that are pullbacks of Specht modules of 
the symmetric group or of the Hecke algebra and therefore obey the usual branching rule. 
For both of them, the restriction contains the 
trivial 
and the sign module for $ KS_2 $ corresponding to the third and fourth 
Specht modules for $ {\cal E}_2(u) $ in the above description.
But the first of them moreover contains the trivial module for $ {\cal H}_2(u) $
corresponding to the first Specht module of the classification, 
whereas the second contains the nontrival one-dimensional 
module for $ {\cal H}_2(u) $ corresponding to the fourth module of the classification.
The question is now how to generalize this to higher $ n $. 

\medskip

The paper treated the representation theory of $ {\cal E}_n(u) $ for
$ u $ generic, where one expects $ {\cal E}_n(u) $ to be semisimple, as
observed above for $ n = 2, 3$.
It is therefore natural to ask for a formal proof of semisimplicity
beyond the cases $ n =2 , 3 $. 
If one had an explicit formula for the dimension of $ S(\Lambda) $ it would be natural 
to try to generalize the above proof for $ n =2 , 3 $. 
On the other hand, in view of the nondegeneracy of the form defined in section 5 
and Wenzl's treatment of the Brauer algebra in [W], 
an attractive alternative approach to proving 
semisimplicity of $ {\cal E}_n(u) $ 
would be to look for an analogue of the Jones basic construction in the setting, using 
the embedding $ {\cal E}_n(u) \subset  {\cal E}_{n+1}(u) $.

\medskip
As already mentioned in section 2, it is possible to define a specialized algebra $ {\cal E}_n(u_0) $, 
for example by choosing $ u_0 $ to be an $l$th root of unity. This should be a 
nonsemisimple algebra. 
A natural first step into the representation theory
of this specialized algebra is 
to show that $ {\cal E}_n(u) $ is a cellular algebra in the sense of [GL]. We firmly believe
that this indeed is the case, but also think that a new set of tools would be needed to establish it. 
In this paper, the tensor module was a crucial 
ingredient in our determination of the rank of $ {\cal E}_n(u) $ 
and so for the completeness of the paper we found it most natural to construct the Specht modules 
inside it. 

\medskip
Finally, the tensor module itself raises the question of determining its endomorphism 
algebra $ \End_{{\cal E}_n(u)} ( V^{ \otimes n}) $ and setting up an analogue of Schur-Weyl duality. 
Given the result of the paper, $ \End_{{\cal E}_n(u)} ( V^{ \otimes n}) $ should be an
interesting combination of 
quantum groups and symmetric groups/Hecke algebras.

\end{document}